\documentclass{article}

\usepackage{amsmath,amssymb,amsthm}
\numberwithin{equation}{section}

\usepackage{geometry}

\usepackage{array}

\usepackage[shortlabels]{enumitem}

\usepackage{tikz}

\usepackage{url}

\usepackage{hyperref}

\usepackage{float}

\usepackage{booktabs}

\usepackage{graphicx}   

\usepackage{indentfirst}

\usepackage{enumitem}
\setlist[itemize]{leftmargin=*} 

\DeclareMathOperator{\supp}{supp}

\DeclareMathOperator{\wt}{wt}

\newtheorem{theorem}{Theorem}[section]

\newtheorem{lemma}[theorem]{Lemma}

\newtheorem{remark}[theorem]{Remark}

\title{On \(q\)-ary Minimal Linear Codes from Hierarchical Posets}\author{
Wen Ma\\
College of Science, Nanjing Forestry University\\
Nanjing 210037, China\\
\texttt{mawen95@126.com}
}

\date{2026}

\begin{document}

\maketitle

\begin{abstract}
In this paper, we study several families of \(q\)-ary linear codes constructed from two-level hierarchical posets over finite fields \(\mathbb{F}_q\), where \(q>2\). By using generating functions and additive characters, we determine their complete weight distributions, dimensions, and minimum distances. This paper focuses on constructing minimal codes that violate the Ashikhmin-Barg condition. By analyzing the weight distributions and support sets of the codewords, we show that a large class of the constructed codes are minimal and  strictly violate the Ashikhmin-Barg condition.

\vspace{1em}
\noindent \textbf{Keywords:} Linear codes;  minimal codes; hierarchical posets; weight distribution; Ashikhmin-Barg condition.
\end{abstract}

\section{Introduction}

Minimal linear codes have attracted considerable attention because of their applications in secret sharing schemes and related areas of information security [Massey1993, Shamir1979, Ding2003]. A fundamental sufficient condition for the minimality of a linear code was established by Ashikhmin and Barg [1998]. In particular, a linear code over $\mathbb F_q$ is minimal if its minimum and maximum nonzero weights, denoted by $w_{\min}$ and $w_{\max}$, satisfy
\[
\frac{w_{\min}}{w_{\max}}>\frac{q-1}{q}.
\]
This condition, commonly referred to as the Ashikhmin--Barg (AB) condition, is sufficient but not necessary. Consequently, the construction of minimal linear codes that do not satisfy this condition has become an interesting problem in coding theory. Various constructions based on Boolean functions and other algebraic methods have been investigated; see, for example, [Ding2016, Ding2018, Heng2018, Lu2021].

Posets provide a natural framework for studying linear codes with additional combinatorial structures. In particular, hierarchical posets have been used to construct linear codes with desirable parameters and weight distributions. Hyun, Kim, Wu and Yue [Hyun2023] studied optimal minimal linear codes arising from posets. Lu, Wu and Cao [lu2021poset] further investigated minimal linear codes constructed from posets. More recently, related constructions have also been considered in [Wu2025, Pelen2026]. These studies demonstrate the usefulness of poset structures in the construction and analysis of minimal linear codes. Nevertheless, the weight distributions and minimality properties of codes associated with collections of order ideals over general finite fields deserve further investigation.

In this paper, we consider two-level hierarchical posets $\mathbb H(m,n)$ and study four classes of $q$-ary linear codes arising from order ideals. The four classes are obtained through two different constructions, each considered in two settings. The first construction is based on the complement $D$ of a union of sets determined by order ideals, giving the codes $C_D$ associated with a single order ideal and with two incomparable order ideals, respectively. The second construction uses a function $f$ defined by the same order-ideal structure, giving the corresponding function-based codes $\mathcal C_f$ in these two settings. Thus, we distinguish throughout the paper between the complement-set construction and the function-based construction.

Our main objective is to determine the weight distributions and investigate the minimality properties of these four classes of linear codes. First, we derive explicit weight distributions and determine the dimensions and parameters of the codes arising from a single order ideal. We then extend the analysis to two incomparable order ideals, where the intersections and differences of the corresponding sets lead to additional cases in the weight distributions. These results provide a systematic description of the four constructions over arbitrary finite fields.

We further investigate the relationship between minimality and the AB condition. For the complement-set construction, we establish conditions under which the resulting codes satisfy the AB condition and hence are minimal. For the function-based construction, we obtain explicit families of minimal linear codes that violate the AB condition. In particular, we give constructions associated with both a single order ideal and two incomparable order ideals. These examples show that the function-based construction can produce minimal linear codes beyond the range covered by the AB sufficient criterion.

The results extend the study of minimal linear codes associated with hierarchical posets to general finite fields. They also illustrate how the choice of defining sets and functions, together with the number of order ideals involved, influences the weight distributions and minimality properties of the resulting codes.

\section{Preliminaries}

We say that \(\mathbb{P} = ([n],\preceq)\) is a partially ordered set (abbreviated as a poset) if \(\preceq\) is a partial order relation on \([n]\), that is, for all \(i,j,k\in [n]\) we have that: (i) \(i\preceq i\); (ii) \(i\preceq j\) and \(j\preceq i\) imply \(i = j\); (iii) \(i\preceq j\) and \(j\preceq k\) imply \(i\preceq k\).

Let \(\mathbb{P} = ([n],\preceq)\) be a poset. Two distinct elements \(i\) and \(j\) in \([n]\) are called comparable if either \(i\preceq j\) or \(j\preceq i\), and incomparable otherwise. It is said that a poset \(\mathbb{P}\) is an anti-chain if every pair of distinct elements is incomparable.

A subset \(I\) of \(\mathbb{P}\), possibly empty, is called an order ideal if \(j\in I\) and \(i\preceq j\) imply \(i\in I\). For a subset \(E\) of \(\mathbb{P}\), the smallest order ideal of \(\mathbb{P}\) containing \(E\) is denoted by \(\langle E\rangle\). For an order ideal \(I\) of \(\mathbb{P}\), we use \(I(\mathbb{P})\) to denote the set of order ideals of \(\mathbb{P}\) that are contained in \(I\). Let \(\mathcal{I} = \{I_1,\ldots ,I_m\}\) be a subset of \(\mathcal{O}_{\mathbb{P}}\). We define

\begin{equation}
\mathcal{I}(\mathbb{P}) = \{J\in \mathcal{O}_{\mathbb{P}}: \text{there exists } I\in\mathcal{I} \text{ such that } J\subseteq I\} = \bigcup_{i = 1}^{m}I_{i}(\mathbb{P}). 
\end{equation}
Then \(\mathcal{I}(\mathbb{P})\) is an order ideal of \(\mathcal{O}_{\mathbb{P}}\) with partial order \(\subseteq\).

Let \(C\) be a linear code over \(\mathbb{F}_q\). For two codewords \(u,v\in C\), we say that \(u\) covers \(v\) if \(\supp(v)\subseteq\supp(u)\). A nonzero codeword \(u\in C\) is called minimal if, for every \(v\in C\) satisfying \(\supp(v)\subseteq\supp(u)\), there exists \(\lambda\in\mathbb{F}_q\) such that \(v=\lambda u\). A linear code \(C\) is called minimal if every nonzero codeword of \(C\) is minimal.

Ashikhmin and Barg \cite{Ashikhmin1998} established the following practical criterion for determining whether a linear code is minimal.

\begin{lemma}
Let \(C\) be a linear code over \(\mathbb{F}_q\), and let \(w_{\min}\) and \(w_{\max}\) denote its minimum and maximum nonzero Hamming weights, respectively. If \(\frac{w_{\min}}{w_{\max}}>\frac{q-1}{q}\), then \(C\) is minimal.
\end{lemma}

The following lemma, which can be derived from the definition of minimal codewords \cite{lu2021}, provides a necessary and sufficient condition for a linear code to be minimal.

\begin{lemma}\label{lem:minimality}
Let \(C\) be a linear code over \(\mathbb{F}_q\). Then \(C\) is a minimal code if and only if, for any two linearly independent codewords \(a,b\in C\),
\[
\sum_{\lambda\in\mathbb{F}_q^*}
\wt(a+\lambda b)
\neq (q-1)\wt(a)-\wt(b).
\]
\end{lemma}

\section{Generating functions for order ideals of \(\mathcal{O}_{\mathbb{P}}\)}

Let \(\mathbb{P}=([n],\preceq)\) be a poset, and let \(\mathcal{I}\subseteq\mathcal{O}_{\mathbb{P}}\). Set
\[
\Delta_{\mathcal{I}}=\{u\in\mathbb{F}_q^n: \supp(u)\in\mathcal{I}(\mathbb{P})\}.
\]
When \(\mathcal{I}=\{I\}\), we write \(\Delta_I\) instead of \(\Delta_{\mathcal{I}}\). For any \(X\subseteq\mathbb{F}_q^n\), define
\[
\mathcal{H}_X(x_1,\ldots,x_n) = \sum_{u\in X}\prod_{i\in\supp(u)}x_i \in\mathbb{Z}[x_1,\ldots,x_n],
\]
where the empty product is \(1\). In particular, \(\mathcal{H}_{\emptyset}=0\) and \(\mathcal{H}_{\{0\}}=1\). If \(X^c=\mathbb{F}_q^n\setminus X\), then
\[
\mathcal{H}_X(x_1,\ldots,x_n)+\mathcal{H}_{X^c}(x_1,\ldots,x_n) = \prod_{i=1}^n(1+(q-1)x_i).
\]

\begin{theorem}\label{rong}
Let \(\mathbb{P}=([n],\preceq)\) be a poset and let \(\mathcal{I}=\{I_1,\ldots,I_k\}\subseteq\mathcal{O}_{\mathbb{P}}\). Then
\begin{equation*}
\mathcal{H}_{\Delta_{\mathcal{I}}}(x_1,\ldots,x_n) = \sum_{\emptyset\neq S\subseteq[k]}(-1)^{|S|+1} \mathcal{H}_{\Delta_{\cap_{j\in S}I_j}}(x_1,\ldots,x_n). 
\end{equation*}
In particular, \(
|\Delta_{\mathcal{I}}| = \sum_{\emptyset\neq S\subseteq[k]}(-1)^{|S|+1} |\Delta_{\cap_{j\in S}I_j}|.
\)
\end{theorem}

\begin{proof}

By definition, \(\Delta_{\mathcal{I}}=\bigcup_{j=1}^k\Delta_{I_j}\). Moreover, for every nonempty \(S\subseteq[k]\),
\[
\bigcap_{j\in S}\Delta_{I_j} = \Delta_{\cap_{j\in S}I_j}.
\]
Let \(\mathbf{1}_X\) denote the indicator function of a subset \(X\subseteq\mathbb{F}_q^n\). The inclusion--exclusion principle gives
\[
\mathbf{1}_{\Delta_{\mathcal{I}}}(u) = \sum_{\emptyset\neq S\subseteq[k]}(-1)^{|S|+1} \mathbf{1}_{\Delta_{\cap_{j\in S}I_j}}(u).
\]
Consequently,
\begin{align*}
\mathcal{H}_{\Delta_{\mathcal{I}}}(x_1,\ldots,x_n)
&=\sum_{u\in\mathbb{F}_q^n} \mathbf{1}_{\Delta_{\mathcal{I}}}(u) \prod_{i\in\supp(u)}x_i\\
&=\sum_{\emptyset\neq S\subseteq[k]}(-1)^{|S|+1} \sum_{u\in\Delta_{\cap_{j\in S}I_j}} \prod_{i\in\supp(u)}x_i\\
&=\sum_{\emptyset\neq S\subseteq[k]}(-1)^{|S|+1} \mathcal{H}_{\Delta_{\cap_{j\in S}I_j}}(x_1,\ldots,x_n).
\end{align*}
Setting \(x_1=\cdots=x_n=1\) gives the cardinality formula.
\end{proof}

\section{Hierarchical posets with two levels}

Let \(m\) and \(n\) be positive integers with \(m \le n\). We say that \(\mathbb{H}(m,n) = ([n], \preceq)\) is a hierarchical poset with two levels if \([n]\) is the disjoint union of two antichains \(U = \{1, \dots, m\}\) and \(V = \{m+1, \dots, n\}\), and \(i \prec j\) whenever \(i \in U\) and \(j \in V\). Its Hasse diagram is shown below. When \(m=n\), we have \(V=\emptyset\), and
\(\mathbb{H}(m,m)\) is an anti-chain.

\begin{center}

\begin{tikzpicture}[scale=1.0]
    \def\dx{2.5}        
    \def\dy{2.5}        
    \def\topoffset{-0.8} 
    \foreach \i/\lbl in {1/1, 2/2, 3/{m-1}, 4/{m}} {
        \node[circle, fill, inner sep=1.5pt, label=below:\(\lbl\)] (B\i) at (\i*\dx, 0) {};
    }
    \foreach \i/\lbl in {1/{m+1}, 2/{m+2}, 3/{n-1}, 4/{n}} {
        \node[circle, fill, inner sep=1.5pt, label=above:\(\lbl\)] (T\i) at (\i*\dx + \topoffset, \dy) {};
    }
    \node at (2.5*\dx - 0.2, 0) {\(\cdots\)};
    \node at (2.5*\dx + \topoffset - 0.2, \dy) {\(\cdots\)};
    \foreach \i in {1,2,3,4} {
        \foreach \j in {1,2,3,4} {
            \draw (B\i) -- (T\j);
        }
    }
\end{tikzpicture}

\end{center}

Every order ideal of \(\mathbb{H}(m,n)\) can be expressed by \(A \cup B\) for \(A \subseteq [m]\), \(B \subseteq [n] \setminus [m]\), and one of the following holds: (i) \(B = \emptyset\); (ii) \(B \neq \emptyset\) and \(A = [m]\).

\begin{lemma}\label{length}

Let \(I = A \cup B\) for \(A \subseteq [m]\), \(B \subseteq [n] \setminus [m]\) be an order ideal of \(\mathbb{P} = \mathbb{H}(m,n)\).

\begin{enumerate}

    \item[(1)] If \(B = \emptyset\), then
 \[
    \mathcal{H}_{\Delta_I}(x_1,\ldots,x_n) = \prod_{i \in A} (1 +(q-1) x_i).
 \]
In particular, we have that \(|\Delta_{\mathcal{I}}| = q^{|A|}\).

    \item[(2)] If \(B \neq \emptyset\), then
\[
    \mathcal{H}_{\Delta_{I}}(x_1, x_2, \dots, x_n) = \prod_{i \in [m]} (1 + (q-1)x_i) +(q-1)^m \prod_{i \in [m]} x_i \left( \prod_{j \in B} (1 + (q-1)x_j) - 1 \right).
    \]
In particular, we have that \(|\Delta_{I}| = q^m + (q-1)^m(q^{|B|} - 1)\).
\end{enumerate}
\end{lemma}

\begin{proof}
The proof is straightforward.
\end{proof}

\section{Weight distributions}

\subsection{Linear codes from hierarchical posets with two levels}
Let \(\mathbb{P}\) be a poset on \([n]\) and let \(
D = \mathbb{F}_q^n\setminus \Delta_{\mathcal{I}}
\) where \(\mathcal{I} = \{I_1, \ldots, I_k\} \subseteq \mathcal{O}_{\mathbb{P}}\). Set
\begin{equation}
C_D = \{c_{D,u} = (u \cdot x)_{x \in D} : u \in \mathbb{F}_q^n\}, 
\end{equation}
where \(\cdot\) denotes Euclidean inner product of two elements in \(\mathbb{F}_q^n\). Then \(C_D\) is a linear code with length \(|D|=q^n-|\Delta_{\mathcal{I}}|\) and its dimension is at most \(n\). Throughout the following, for a subset \(X \subseteq [n]\), \(\delta(u|X)\) denotes the indicator function satisfying \(\delta(u|X) = 1\) if and only if \(supp(u) \cap X = \emptyset\). We also use \(\delta_{0,u}\) to denote the Kronecker delta, i.e., \(\delta_{0,u} = 1\) when \(u=0\), and \(0\) otherwise.

Let \(\chi\) be a nontrivial additive character of the field \(\mathbb{F}_q\). The Hamming weight of \(c_{D,u}\) is
\begin{equation}\label{eq:weight-CD}
\begin{aligned}
\wt(c_{D,u})
&= |D| - \frac{1}{q}\sum_{y\in\mathbb{F}_q}\sum_{x\in D}\chi(y(u\cdot x))\\
&= \frac{q-1}{q}|D| - \frac{q-1}{q}q^n\delta_{0,u}
   +\frac{1}{q}\sum_{y\in\mathbb{F}_q^*}\sum_{x\in\Delta_{\mathcal I}}\chi(y(u\cdot x))\\
&= \frac{q-1}{q}\left(|D|+\mathcal H_{\Delta_{\mathcal I}}(w_1,\ldots,w_n)-q^n\delta_{0,u}\right).
\end{aligned}
\end{equation}
where
\[
w_i=\frac{1}{q-1}\sum_{t\in\mathbb F_q^*}\chi(u_i t)
=\begin{cases}
1,&u_i=0,\\[1mm]
-\dfrac{1}{q-1},&u_i\ne0.
\end{cases}
\]
The value of the average is independent of the nonzero scalar multiplying \(u_i\), so the same tuple \((w_1,\ldots,w_n)\) occurs for every \(y\in\mathbb F_q^*\). In particular, the displayed formula is valid for both \(u=0\) and \(u\ne0\).

\begin{theorem}\label{T1}

Let \(\mathbb{H}(m,n)\) be a hierarchical poset with two levels and \(I = A \cup B\) an order ideal of \(\mathbb{H}(m,n)\) for \(A \subseteq [m]\), \(B \subseteq [n] \setminus [m]\). Set \(\mathcal{I} = \{I\}\).

\begin{enumerate}

    \item[(1)] If \(B=\emptyset\) and \(A\ne[n]\), then the weight distribution of \(C_D\) is given in Table~\ref{tab:weight1}.
    \item[(2)] If \(B\ne\emptyset\) and it is not the exceptional case \(m=1,\ B=[n]\setminus[m]\), then the weight distribution of \(C_D\) is given in Table~\ref{tab:weight2}. If \(m=1\) and \(B=[n]\setminus[m]\), then \(C_D\) is the \(q\)-ary simplex code with weight distribution \(0\) (frequency \(1\)) and \((q-1)q^{n-2}\) (frequency \(q^{n-1}-1\)).

\end{enumerate}

\end{theorem}

\begin{proof}

Let \(\mathbb{P} = \mathbb{H}(m,n)\).  For any \(u = (u_1, u_2, \ldots, u_n) \in \mathbb{F}_q^n\), we write \(u = (v, w)\), where \(v = (u_1, \ldots, u_m)\) and \(w = (u_{m+1}, \ldots, u_n)\).

(1) By Lemma \ref{length}, the length of the code \(C_D\) is \(q^n - |\Delta_{\mathcal{I}}| = q^n - q^{|A|}\). Note that
\[
\mathcal{H}_{\Delta_{\mathcal{I}}}(x_1, x_2, \dots, x_n) = \mathcal{H}_{\Delta_{A}}(x_1, x_2, \dots, x_n) = \prod_{i \in A} (1 +(q-1) x_i).
\]
Substituting \(w_i = \frac{1}{q-1}\sum_{t\in\mathbb{F}_q^*} \chi(u_i ty)\) into the expression, we obtain
\begin{equation*}
1+(q-1)w_i = 
\begin{cases} 
q, & u_i = 0, \\[2mm]
0, & u_i \neq 0.
\end{cases}
\end{equation*}
Consequently, evaluating the generating function gives
\begin{align*}
\mathcal{H}_{\Delta_{\mathcal{I}}}\left( w_1, \dots, w_n \right)
= \prod_{i \in A} (1 + (q-1)w_i) 
=q^{|A|}\delta(u|A).
\end{align*}
For \(u\ne0\), substituting this result into the weight formula \eqref{eq:weight-CD}, we obtain
\begin{equation*}
\wt(c_{D,u}) = (q-1)q^{n-1} - \frac{q-1}{q}q^{|A|} \bigl(1 - \delta(u|A)\bigr).
\end{equation*}
For \(u=0\), we have \(\wt(c_{D,0})=0\). The result then follows.

(2) Let \(B \neq \emptyset\). If \(m=1\) and \(B=[n]\setminus[m]\), then \(D=\{(0,w):w\in\mathbb F_q^{n-1}\setminus\{0\}\}\), so \(C_D\) is the simplex code. We may assume this exceptional case does not occur. By Lemma \ref{length}, the length of the code \(C_D\) is  
\[
N=q^n - |\Delta_{\mathcal{I}}|= q^n - q^m - (q-1)^m(q^{|B|} - 1).
\]
Note that 
\[
\mathcal{H}_{\Delta_{\mathcal{I}}}(x_1, x_2 \dots, x_n) = \prod_{i \in [m]} (1 + (q-1)x_i) +(q-1)^m \prod_{i \in [m]} x_i \left( \prod_{j \in B} (1 + (q-1)x_j) - 1 \right).
\]
Therefore, 
\begin{align*}
&  \mathcal{H}_{\Delta_{\mathcal{I}}}(w_1, w_2, \dots, w_n) \\
=& \prod_{i \in [m]} (1 + (q-1)w_i) + (q-1)^m\prod_{i \in [m]} w_i \left( \prod_{j \in B} (1 + (q-1)w_j) - 1 \right) \\
=& q^m \delta(u|[m]) + (q-1)^m \left( (-1)^{\wt(v)} (q-1)^{-\wt(v)} \right) \left( q^{|B|}\delta(u|B) - 1 \right) \\
=& q^m \delta(u|[m]) + (-1)^{\wt(v)} (q-1)^{m-\wt(v)} \left( q^{|B|}\delta(u|B) - 1 \right).
\end{align*}
For \(u\ne0\), substituting the above expression into the weight formula \eqref{eq:weight-CD}, we obtain
\begin{align*}
\wt(c_{D,u})
=\frac{q-1}{q}\Big[
&N+q^m\delta(u|[m])+(-1)^{\wt(v)}(q-1)^{m-\wt(v)}
\big(q^{|B|}\delta(u|B)-1\big)
\Big],
\end{align*}
where \(N=q^n-q^m-(q-1)^m(q^{|B|}-1)\).

\begin{enumerate}[label=\textbf{Case \arabic*:}, leftmargin=*]

\item \(u=0\). The codeword is the zero vector, so \(\wt(c_{D,0})=0\).
\item \(v=0\) and \(w|_B=0\), with \(w\neq0\). Here \(\delta(u|[m])=1\) and \(\delta(u|B)=1\). Then
\[
\wt(c_{D,u}) =\frac{q-1}{q} \left[N+q^m+(q-1)^m(q^{|B|}-1)\right] = (q-1)q^{n-1}.
\]
\item \(v=0\) and \(w|_B\neq0\). Here \(\delta(u|[m])=1\) and \(\delta(u|B)=0\). Thus,
\[
\wt(c_{D,u}) =\frac{q-1}{q} \left[N+q^m-(q-1)^m\right] = (q-1)q^{n-1} -(q-1)^{m+1}q^{|B|-1}.
\]
\item \(v\neq0\) and \(w|_B=0\). Here \(\delta(u|[m])=0\) and \(\delta(u|B)=1\). Let \(r=\wt(v)\geq1\). Then
\[
\wt(c_{D,u}) = \frac{q-1}{q} \left[ N+(-1)^r(q-1)^{m-r}(q^{|B|}-1) \right].
\]
\item \(v\neq0\) and \(w|_B\neq0\). Here \(\delta(u|[m])=0\) and \(\delta(u|B)=0\). Let \(r=\wt(v)\geq1\). Then
\[
\wt(c_{D,u}) = \frac{q-1}{q} \left[ N-(-1)^r(q-1)^{m-r} \right].
\]

\end{enumerate}

\begin{table}[htbp]

\centering

\caption{Weight distribution for Theorem \ref{T1} (1)}

\label{tab:weight1}

\begin{tabular}{cc}
\toprule
Weight & Frequency \\
\midrule
\(0\) & \(1\) \\
\((q-1)q^{n-1}\) & \(q^{n-|A|} - 1\) \\
\(\dfrac{q-1}{q}(q^n-q^{|A|})\) & \(q^n - q^{n-|A|}\) \\
\bottomrule
\end{tabular}

\end{table}

\begin{table}[htbp]
\centering
\renewcommand{\arraystretch}{1.2}
\caption{Weight distribution for Theorem \ref{T1} (2)}
\label{tab:weight2}
\small
\begin{tabular}{cc}
\toprule
Weight &Frequency \\ 
\midrule
\(0\) & \(1\) \\
\((q-1)q^{n-1}\)
& \(q^{n-m-|B|}-1\) \\
\((q-1)q^{n-1}-(q-1)^{m+1}q^{|B|-1}\)
& \((q^{|B|}-1)q^{n-m-|B|}\) \\
\(\displaystyle
\frac{q-1}{q}
\left[
q^n-q^m-
\left((q-1)^m-(-1)^r(q-1)^{m-r}\right)
(q^{|B|}-1)
\right]
\)
&
\(\displaystyle
\binom{m}{r}(q-1)^r q^{n-m-|B|}
\) \\
\(\displaystyle
\frac{q-1}{q}
\left[
q^n-q^m-(q-1)^m(q^{|B|}-1)
-(-1)^r(q-1)^{m-r}
\right]
\)
&
\(\displaystyle
\binom{m}{r}(q-1)^r
(q^{|B|}-1)q^{n-m-|B|}
\) \\
\bottomrule
\end{tabular}%
\begin{flushleft}
\footnotesize
Here \(r=\wt(v)\) ranges over \(1\leq r\leq m\). Each integer \(r\) yields a specific weight. If different values of \(r\) produce the same weight, their corresponding frequencies should be added.
\end{flushleft}
\end{table}
\end{proof}

\begin{remark}
The parameters and dimensions of the codes in Theorem~\ref{T1} are as follows.

\begin{enumerate}[(1)]
\item If \(B=\emptyset\) and \(A=\emptyset\), then \(C_D\) is the \(q\)-ary simplex code with parameters
\(
[q^n-1,n,(q-1)q^{n-1}].
\)
If \(1\le |A|\le n-1\), then
\(
C_D\text{ has parameters }[q^n-q^{|A|},n,(q-1)(q^{n-1}-q^{|A|-1})].
\)
These codes attain the Griesmer bound. If \(A=[n]\), then \(D=\emptyset\) and \(C_D\) is the zero code.

\item If \(B\ne\emptyset\), put \(b=|B|\) and \(c=n-m-b\). The dimension is \(n\), except when \(m=1\) and \(c=0\), in which case \(D=\{(0,w):w\in\mathbb F_q^{n-1}\setminus\{0\}\}\) and \(C_D\) is the simplex code with parameters
\(
[q^{n-1}-1,n-1,(q-1)q^{n-2}].
\)
\end{enumerate}
\end{remark}

Consider two distinct order ideals \(I_1\) and \(I_2\). If \(I_1\) contains only lower-layer elements while \(I_2\) contains upper-layer elements, then the hierarchical structure forces \(I_1 \subset I_2\). In this case, \(\mathcal{I}(\mathbb{P})\) reduces to \(I_2(\mathbb{P})\), making the construction equivalent to the single-ideal case already handled in Theorem \ref{T1} (2). Therefore, our subsequent analysis will focus solely on the cases where neither ideal is contained in the other.

\begin{theorem}\label{T2}

Consider a two-level hierarchical poset \(\mathbb{H}(m,n)\).  Let \(I_1 = A_1 \cup B_1\) and \(I_2 = A_2 \cup B_2\) be two distinct order ideals of \(\mathbb{H}(m,n)\), where \(A_i \subseteq [m]\), \(B_i \subseteq [n] \setminus [m]\), \(i = 1, 2\) and \(I_1 \nsubseteq I_2\), \(I_2 \nsubseteq I_1\). Set \(\mathcal{I} = \{I_1, I_2\}\).

\begin{enumerate}

    \item[(1)] If \(B_1=B_2 = \emptyset\), then the weight distribution of the code \(C_D\) is given in Table~\ref{tab:weight3_q}.

\item[(2)] If $B_1,B_2\ne\emptyset$, then the weight distribution of the code $C_D$ is given in Table~\ref{tab:weight4_q}.
\end{enumerate}
\end{theorem}

\begin{proof}
(1) Suppose \(B_1=B_2=\emptyset\). Since \(I_1\cap I_2=A_1\cap A_2\), inclusion--exclusion gives
\[
|D|=q^n-q^{|A_1|}-q^{|A_2|}+q^{|A_1\cap A_2|}=N.
\]
Moreover,
\begin{align*}
\mathcal{H}_{\Delta_{\mathcal{I}}}(w_1,\ldots,w_n)
={}&q^{|A_1|}\delta(u|A_1)+q^{|A_2|}\delta(u|A_2)-q^{|A_1\cap A_2|}\delta(u|A_1\cap A_2).
\end{align*}
We distinguish the following cases:
\begin{enumerate}[label=\textbf{Case \arabic*:}, leftmargin=*]
\item If \(u=0\), then \(\wt(c_{D,u})=0\); the frequency is \(1\).
\item If \(u\ne0\) and \(\delta(u|A_1)=\delta(u|A_2)=1\), then
\(
\wt(c_{D,u})=(q-1)q^{n-1}.
\)
\item If \(\delta(u|A_1)=1\) and \(\delta(u|A_2)=0\), then \(\delta(u|A_1\cap A_2)=1\), and
\[
\wt(c_{D,u})=\frac{q-1}{q}(N+q^{|A_1|}-q^{|A_1\cap A_2|})
=\frac{q-1}{q}(q^n-q^{|A_2|}).
\]
\item If \(\delta(u|A_1)=0\) and \(\delta(u|A_2)=1\), then \(\delta(u|A_1\cap A_2)=1\), and
\[
\wt(c_{D,u})=\frac{q-1}{q}(q^n-q^{|A_1|}).
\]
\item If \(\delta(u|A_1)=\delta(u|A_2)=0\) and \(\delta(u|A_1\cap A_2)=1\), then
\[
\wt(c_{D,u})=\frac{q-1}{q}(N-q^{|A_1\cap A_2|})
=\frac{q-1}{q}(q^n-q^{|A_1|}-q^{|A_2|}).
\]
\item If \(\delta(u|A_1\cap A_2)=0\), then \(\delta(u|A_1)=\delta(u|A_2)=0\), and
\[
\wt(c_{D,u})=\frac{q-1}{q}N.
\]
\end{enumerate}

\begin{table}[htbp]

\centering

\renewcommand{\arraystretch}{1.6}

\caption{Weight distribution in Theorem~\ref{T2}(1) over \(\mathbb{F}_q\)}

\label{tab:weight3_q}

\small 

\begin{tabular}{cc}
\toprule
\textbf{Weight} & \textbf{Frequency}\\
\midrule
\(0\) & \(1\)\\
\((q-1)q^{n-1}\) & \(q^{n-|A_1\cup A_2|}-1\)\\
\(\dfrac{q-1}{q}(q^n-q^{|A_2|})\)
& \(q^{n-|A_1\cup A_2|}(q^{|A_2\setminus A_1|}-1)\)\\
\(\dfrac{q-1}{q}(q^n-q^{|A_1|})\)
& \(q^{n-|A_1\cup A_2|}(q^{|A_1\setminus A_2|}-1)\)\\
\(\dfrac{q-1}{q}(q^n-q^{|A_1|}-q^{|A_2|})\)
& \(q^{n-|A_1\cup A_2|}(q^{|A_1\setminus A_2|}-1)(q^{|A_2\setminus A_1|}-1)\)\\
\(\dfrac{q-1}{q}(q^n-q^{|A_1|}-q^{|A_2|}+q^{|A_1\cap A_2|})\)
& \((q^{|A_1\cap A_2|}-1)q^{n-|A_1\cap A_2|}\)\\
\bottomrule
\end{tabular}

\end{table}
(2) Suppose that \(B_1\neq\emptyset\) and \(B_2\neq\emptyset\). By Lemma \ref{rong}, we have
\begin{align*}
|\Delta_{\mathcal{I}}|
=|\Delta_{I_1}|+|\Delta_{I_2}|-|\Delta_{I_1\cap I_2}|
=q^m+(q-1)^m
\bigl(q^{|B_1|}+q^{|B_2|}
-q^{|B_1\cap B_2|}-1\bigr).
\end{align*}
Therefore, the length of \(\mathcal{C}_D\) is
\[
\begin{aligned}
|D|=q^n-q^m-(q-1)^m
\bigl(q^{|B_1|}+q^{|B_2|}
-q^{|B_1\cap B_2|}-1\bigr).
\end{aligned}
\]
Moreover,
\begin{align*}
& \mathcal{H}_{\Delta_{\mathcal{I}}}(w_1, w_2, \dots, w_n) \\\ 
={}& \prod_{i \in [m]} (1 + (q-1)w_i) + (q-1)^m\prod_{i \in [m]} w_i \Bigg( \prod_{j \in B_1} (1 + (q-1)w_j) \\
& + \prod_{j \in B_2} (1 + (q-1)w_j) - \prod_{j \in B_1\cap B_2} (1 + (q-1)w_j) - 1 \Bigg) \\[2mm]
={}& q^m \delta(u|[m]) + (-1)^{\wt(v)} (q-1)^{m-\wt(v)} \Bigg( \prod_{j \in B_1} (1 + (q-1)w_j) \\
& + \prod_{j \in B_2} (1 + (q-1)w_j) - \prod_{j \in B_1\cap B_2} (1 + (q-1)w_j) - 1 \Bigg) \\[2mm]
={}& q^m \delta(u|[m]) + (-1)^{\wt(v)} (q-1)^{m-\wt(v)} \Big( q^{|B_1|}\delta(u|B_1) + q^{|B_2|}\delta(u|B_2) \\
& - q^{|B_1\cap B_2|}\delta(u|B_1\cap B_2) - 1 \Big)
\end{align*}
For \(u=0\), we have \(\wt(c_{D,0})=0\).  For convenience, we define
\[
\Delta_B(u)=
q^{|B_1|}\delta(u|B_1)
+q^{|B_2|}\delta(u|B_2)-q^{|B_1\cap B_2|}
\delta(u|B_1\cap B_2)-1.
\]

\medskip

\noindent \textbf{Part 1: \(v=0\).} For \(u\neq0\), equation~\eqref{eq:weight-CD} becomes
\[
\wt(c_{D,u})=\frac{q-1}{q}\big[N+q^m+(q-1)^m\Delta_B\big].
\]

\begin{enumerate}[label=\textbf{Case \arabic*:}, leftmargin=*]

\item \(u\neq0\) and \(\delta(u|B_1)=\delta(u|B_2)=1\). Then \(\Delta_B(u) = q^{|B_1|} + q^{|B_2|} - q^{|B_1\cap B_2|} - 1\), and \[
\wt(c_{D,u})=(q-1)q^{n-1}.
\]

\item \(\delta(u|B_1)=1\) and \(\delta(u|B_2)=0\). Then \(\delta(u|B_1\cap B_2)=1\), \(\Delta_B(u)=q^{|B_1|}-q^{|B_1\cap B_2|}-1\) and 
\[
\wt(c_{D,u})=\frac{q-1}{q}\bigl(q^n-(q-1)^m q^{|B_2|}\bigr).
\]

\item \(\delta(u|B_1)=0\) and \(\delta(u|B_2)=1\). Then \(\delta(u|B_1\cap B_2)=1\) and \(\Delta_B(u)=q^{|B_2|}-q^{|B_1\cap B_2|}-1\). Hence
\[
\wt(c_{D,u})=\frac{q-1}{q}\bigl(q^n-(q-1)^m q^{|B_1|}\bigr).
\]

\item \(\delta(u|B_1)=\delta(u|B_2)=0\) and \(\delta(u|B_1\cap B_2)=1\). Then \(\Delta_B(u)=-q^{|B_1\cap B_2|}-1\), and
\[
\wt(c_{D,u})=\frac{q-1}{q}
\bigl(q^n-(q-1)^m(q^{|B_1|}+q^{|B_2|})\bigr).
\]

\item \(\delta(u|B_1\cap B_2)=0\). Then \(\delta(u|B_1)=\delta(u|B_2)=0\) and \(\Delta_B(u)=-1\). Thus
\[
\wt(c_{D,u})=\frac{q-1}{q}
\bigl(q^n-(q-1)^m(q^{|B_1|}+q^{|B_2|}-q^{|B_1\cap B_2|})\bigr).
\]

\end{enumerate}

\medskip

\noindent \textbf{Part 2: \(v\neq0\).}  Fix \(r=\wt(v)\) with \(1\le r\le m\). There are \(\binom{m}{r}(q-1)^r\) choices for \(v\). Since \(v\neq0\), no additional subtraction of \(1\) is needed to exclude the zero vector \(u\). Equation~\eqref{eq:weight-CD} becomes
\[
\wt(c_{D,u})=\frac{q-1}{q}
\big[N+(-1)^r(q-1)^{m-r}\Delta_B(u)\big].
\]

\begin{enumerate}[label=\textbf{Case \arabic*:}, leftmargin=*, start=6]

\item \(\delta(u|B_1)=\delta(u|B_2)=1\).
Then \(\Delta_B (u)= q^{|B_1|} + q^{|B_2|} - q^{|B_1\cap B_2|} - 1\), and

\[
\wt(c_{D,u})=\frac{q-1}{q}
\big[N+(-1)^r(q-1)^{m-r}\big(q^{|B_1|} + q^{|B_2|} - q^{|B_1\cap B_2|} - 1\big)\big].
\]

\item \(\delta(u|B_1)=1\) and \(\delta(u|B_2)=0\). Then \(\Delta_B(u)=q^{|B_1|}-q^{|B_1\cap B_2|}-1\) and
\[
\wt(c_{D,u})=\frac{q-1}{q}
\big[N+(-1)^r(q-1)^{m-r}(q^{|B_1|}-q^{|B_1\cap B_2|}-1)\big].
\]

\item \(\delta(u|B_1)=0\) and \(\delta(u|B_2)=1\). Then \(\Delta_B(u)=q^{|B_2|}-q^{|B_1\cap B_2|}-1\) and
\[
\wt(c_{D,u})=\frac{q-1}{q}
\big[N+(-1)^r(q-1)^{m-r}(q^{|B_2|}-q^{|B_1\cap B_2|}-1)\big].
\]

\item \(\delta(u|B_1)=\delta(u|B_2)=0\) and \(\delta(u|B_1\cap B_2)=1\). Then \(\Delta_B(u)=-q^{|B_1\cap B_2|}-1\) and
\[
\wt(c_{D,u})=\frac{q-1}{q}
\big[N-(-1)^r(q-1)^{m-r}(q^{|B_1\cap B_2|}+1)\big].
\]

\item \(\delta(u|B_1\cap B_2)=0\). Then \(\Delta_B(u)=-1\) and
\[
\wt(c_{D,u})=\frac{q-1}{q}
\big[N-(-1)^r(q-1)^{m-r}\big].
\]
\end{enumerate}

\begin{table}[htbp]

\centering

\caption{Weights and frequencies for Theorem~\ref{T2}(2)}

\label{tab:weight4_q}

\renewcommand{\arraystretch}{2.0} 

\resizebox{\textwidth}{!}{%

\begin{tabular}{cc}
\toprule
\textbf{Weight} & \textbf{Frequency} \\
\midrule
\(0\) & \(1\) \\
\((q-1)q^{n-1}\) & \(q^{n-m-|B_1 \cup B_2|}-1\) \\
\(\dfrac{q-1}{q}\big(q^n-(q-1)^m q^{|B_2|}\big)\) & \(q^{n-m-|B_1 \cup B_2|}(q^{|B_2 \setminus B_1|}-1)\) \\
\(\dfrac{q-1}{q}\big(q^n-(q-1)^m q^{|B_1|}\big)\) & \(q^{n-m-|B_1 \cup B_2|}(q^{|B_1 \setminus B_2|}-1)\) \\
\(\dfrac{q-1}{q}\big(q^n-(q-1)^m(q^{|B_1|}+q^{|B_2|})\big)\) & \(q^{n-m-|B_1 \cup B_2|}(q^{|B_1 \setminus B_2|}-1)(q^{|B_2 \setminus B_1|}-1)\) \\
\(\dfrac{q-1}{q}\big(q^n-(q-1)^m(q^{|B_1|}+q^{|B_2|}-q^{|B_1\cap B_2|})\big)\) & \((q^{|B_1\cap B_2|}-1)q^{n-m-|B_1\cap B_2|}\) \\
\(\dfrac{q-1}{q}\Big[ q^n - q^m - \big( (q-1)^m - (-1)^r(q-1)^{m-r} \big)(q^{|B_1|} + q^{|B_2|} - q^{|B_1\cap B_2|} - 1) \Big]\) & \(\binom{m}{r}(q-1)^r q^{n-m-|B_1 \cup B_2|}\) \\
\(\dfrac{q-1}{q}\Big[ q^n - q^m - (q-1)^m(q^{|B_1|}+q^{|B_2|}-q^{|B_1\cap B_2|}-1) + (-1)^r(q-1)^{m-r}(q^{|B_1|}-q^{|B_1\cap B_2|}-1) \Big]\) & \(\binom{m}{r}(q-1)^r q^{n-m-|B_1 \cup B_2|}(q^{|B_2 \setminus B_1|}-1)\) \\
\(\dfrac{q-1}{q}\Big[ q^n - q^m - (q-1)^m(q^{|B_1|}+q^{|B_2|}-q^{|B_1\cap B_2|}-1) + (-1)^r(q-1)^{m-r}(q^{|B_2|}-q^{|B_1\cap B_2|}-1) \Big]\) & \(\binom{m}{r}(q-1)^r q^{n-m-|B_1 \cup B_2|}(q^{|B_1 \setminus B_2|}-1)\) \\
\(\dfrac{q-1}{q}\Big[ q^n - q^m - (q-1)^m(q^{|B_1|}+q^{|B_2|}-q^{|B_1\cap B_2|}-1) - (-1)^r(q-1)^{m-r}(q^{|B_1\cap B_2|}+1) \Big]\) & \(\binom{m}{r}(q-1)^r q^{n-m-|B_1 \cup B_2|}(q^{|B_1 \setminus B_2|}-1)(q^{|B_2 \setminus B_1|}-1)\) \\
\(\dfrac{q-1}{q}\Big[ q^n - q^m - (q-1)^m(q^{|B_1|}+q^{|B_2|}-q^{|B_1\cap B_2|}-1) - (-1)^r(q-1)^{m-r} \Big]\) & \(\binom{m}{r}(q-1)^r(q^{|B_1\cap B_2|}-1)q^{n-m-|B_1\cap B_2|}\) \\
\bottomrule
\end{tabular}%

}

\vspace{2mm}

\begin{flushleft}
\footnotesize
Here \(r=\wt(v)\) with \(1\le r\le m\). If two or more entries yield the same weight, their corresponding frequencies are added together.
\end{flushleft}
\end{table}
\end{proof}

 \begin{remark}
For the codes in Theorem~\ref{T2}, the dimension is
\[
\dim(C_D)=n.
\]
Indeed, in both cases we have
\(
\operatorname{span}_{\mathbb F_q}(D)=\mathbb F_q^n.
\)
Consequently, the evaluation map
\[
u\longmapsto c_{D,u}
\]
is injective.
\end{remark}

\subsection{Linear codes from \(q\)-ary functions}

By a \(q\)-ary function, we mean a function from \(\mathbb{F}_q^n\) to \(\mathbb{F}_q\). Let \(f:\mathbb{F}_q^n\to\mathbb{F}_q\) be a function satisfying \(f(0)=0\) and \(f(x)\neq0\) for at least one \(x\in\mathbb{F}_q^n\). We define a linear code associated with \(f\) by
\begin{equation}
\mathcal{C}_f = \left\{ c_f(s,u) = (sf(x) + u \cdot x)_{x\in\mathbb{F}_q^n\setminus\{0\}} : s\in\mathbb{F}_q,\, u\in\mathbb{F}_q^n \right\}. 
\end{equation}
The length of \(\mathcal{C}_f\) is \(q^n-1\), and its dimension is at most \(n+1\).

Let \(\mathbb{P}\) be a poset on \([n]\), and let \(\mathcal{I}=\{I_1,\ldots,I_k\}\subseteq\mathcal{O}_{\mathbb{P}}\). Consider the indicator function \(f:\mathbb{F}_q^n\to\mathbb{F}_q\) defined by
\[
f(x) = \begin{cases} 1, & x \in \Delta_{\mathcal{I}} \setminus \{0\}, \\\ 0, & \text{otherwise}. \end{cases}
\]

We assume that \(|\Delta_{\mathcal{I}}|>1\), so that \(f\) is not identically zero. Let \(\chi\) be a nontrivial additive character of \(\mathbb{F}_q\). The generalized Walsh--Hadamard transform of \(f\) with respect to \(\chi\) is defined by
\[
S_f(u)=
\sum_{x\in\mathbb{F}_q^n}
\chi(f(x)+u\cdot x).
\]
Such transforms have been studied in the context of \(q\)-ary plateaued functions by Mesnager et al.~\cite{Mesnager2019}, following earlier work on bent polynomials over finite fields by Coulter and
Matthews~\cite{Coulter1997}.

By the definition of \(f\), we can split the generalized Walsh--Hadamard transform \(S_f(u)\) over the whole space \(\mathbb{F}_q^n\) and the subset \(\Delta_{\mathcal{I}} \setminus \{0\}\):
\begin{align*}
S_f(u) 
&= \sum_{x \in \mathbb{F}_q^n} \chi(f(x) + u \cdot x) \\
&= \sum_{x \in \mathbb{F}_q^n} \chi(u \cdot x) + \sum_{x \in \Delta_{\mathcal{I}} \setminus \{0\}} \bigl( \chi(1 + u \cdot x) - \chi(u \cdot x) \bigr) \\
&= \sum_{x \in \mathbb{F}_q^n} \chi(u \cdot x) + (\chi(1) - 1) \sum_{x \in \Delta_{\mathcal{I}} \setminus \{0\}} \chi(u \cdot x).
\end{align*}

Using the orthogonality of the character \(\chi\), the first sum is \(q^n \delta_{0,u}\). For the second sum, since the generating function contains the zero vector, we have \(\sum_{x \in \Delta_{\mathcal{I}} \setminus \{0\}} \chi(u \cdot x) = \mathcal{H}_{\Delta_{\mathcal{I}}}(w_1, \dots, w_n) - 1\). Thus, 
\[
S_f(u) = q^n\delta_{0,u} + (\chi(1)-1) \bigl( \mathcal{H}_{\Delta_{\mathcal{I}}}(w_1,\ldots,w_n) - 1 \bigr). 
\]
The Hamming weight of the codeword \(c_f(s,u)\) is calculated using the orthogonality of characters:
\begin{align*}
\wt(c_f(s,u)) 
&= \sum_{x \in \mathbb{F}_q^n\setminus\{0\}} \left( 1 - \frac{1}{q} \sum_{y\in\mathbb{F}_q} \chi(y(sf(x)+u \cdot x)) \right) \\
&= (q^n - 1) - \frac{1}{q} \sum_{y\in\mathbb{F}_q} \sum_{x\in\mathbb{F}_q^n\setminus\{0\}} \chi(y(sf(x)+u \cdot x))\\
&= (q^n - 1) - \frac{1}{q} \left( (q^n - 1) + \sum_{y\in\mathbb{F}_q^*} \bigl(S_{ysf}(yu)-1\bigr) \right) \\
&= (q-1)q^{n-1} - \frac{1}{q} \sum_{y\in\mathbb{F}_q^*} S_{ysf}(yu).
\end{align*}
Substituting the formula of \(S_{ysf}(yu)\) into the equation:
\begin{align*}
S_{ysf}(yu) = q^n\delta_{0,u} + (\chi(ys)-1) \left( \mathcal{H}_{\Delta_{\mathcal{I}}}(w_1,\ldots,w_n) - 1 \right),
\end{align*}
we can evaluate the sum as
\begin{align*}
\sum_{y\in\mathbb{F}_q^*} S_{ysf}(yu) 
&= \sum_{y\in\mathbb{F}_q^*} q^n\delta_{0,u} + \left( \mathcal{H}_{\Delta_{\mathcal{I}}}(w_1,\ldots,w_n) - 1 \right) \sum_{y\in\mathbb{F}_q^*} (\chi(ys)-1) \\
&= (q-1)q^n\delta_{0,u} + \left( \mathcal{H}_{\Delta_{\mathcal{I}}}(w_1,\ldots,w_n) - 1 \right) \sum_{y\in\mathbb{F}_q^*} (\chi(ys)-1).
\end{align*}
Since the additive character \(\chi\) sums to \(-1\) over non-zero elements, we have
\[
\sum_{y\in\mathbb{F}_q^*} (\chi(ys)-1) = \begin{cases} 0, & s=0, \\\ -q, & s\neq0. \end{cases}
\]
Thus, substituting back yields
\[
\wt(c_f(s,u)) = (q-1)q^{n-1} - \frac{1}{q} \Bigl[ (q-1)q^n\delta_{0,u} - q(1 - \delta_{s,0})\bigl( \mathcal{H}_{\Delta_{\mathcal{I}}}(w_1,\ldots,w_n) - 1 \bigr) \Bigr],
\]
which simplifies to
\begin{equation}\label{eq:weight-Cf}
\wt(c_f(s,u)) = (q-1)q^{n-1}(1 - \delta_{0,u}) + (1 - \delta_{s,0})\bigl( \mathcal{H}_{\Delta_{\mathcal{I}}}(w_1,\ldots,w_n) - 1 \bigr).
\end{equation}

Finally, evaluating this formula in all four cases gives

\[
\wt(c_f(s,u)) = \begin{cases} 
0, & s=0,\, u=0, \\[2mm]
(q-1)q^{n-1}, & s=0,\, u\neq0, \\[2mm]
|\Delta_{\mathcal{I}}|-1, & s\neq0,\, u=0, \\[2mm]
(q-1)q^{n-1} + \mathcal{H}_{\Delta_{\mathcal{I}}}(w_1,\ldots,w_n) - 1, & s\neq0,\, u\neq0.
\end{cases}
\]

\begin{theorem}\label{T3}

 \textit{Let \(\mathbb{H}(m,n)\) be a hierarchical poset with two levels and let \(I = A \cup B\) be a
nonempty order ideal of \(\mathbb{H}(m,n)\), where \(A \subseteq [m]\) and \(B \subseteq [n] \setminus [m]\). Set \(\mathcal{I} = \{I\}\). For the \(q\)-ary function \(f: \mathbb{F}_q^n \to \mathbb{F}_q\) defined with respect to \(\mathcal{I}\), the weight distribution of the associated linear code \(\mathcal{C}_f\) is given as follows:}

\begin{enumerate}[(1)]

   \item If $B=\emptyset$, then $A\neq\emptyset$ and the weight distribution of the code
$\mathcal C_f$ is given in Table~\ref{tab:c_f_case1}.

    \item \textit{If \(B \neq \emptyset\), then the weight distribution of the code \(\mathcal{C}_f\) is given in Table \ref{tab:c_f_case2}.}

\end{enumerate}

\end{theorem}

\begin{table}[htbp]

\centering
\caption{Weight distribution in Theorem~\ref{T3}(1): $B=\emptyset$}

\label{tab:c_f_case1}

\renewcommand{\arraystretch}{1.5}

\small
\begin{tabular}{cc}
\toprule
\textbf{Weight} & \textbf{Frequency} \\
\midrule
\(0\) & \(1\) \\
\((q-1)q^{n-1}\) & \(q^n - 1\) \\
\(q^{|A|} - 1\) & \(q - 1\) \\
\((q-1)q^{n-1} + q^{|A|} - 1\) & \((q-1)(q^{n-|A|} - 1)\) \\
\((q-1)q^{n-1} - 1\) & \((q-1)(q^{|A|} - 1)q^{n-|A|}\) \\
\bottomrule
\end{tabular}

\end{table}

\begin{proof}
(1) Suppose that \(B = \emptyset\) and \(A\ne\emptyset\). Then \(|\Delta_{\mathcal{I}}| = q^{|A|}\) and the generating function simplifies to
\[
\mathcal{H}_{\Delta_{\mathcal{I}}}(x_1,\ldots,x_n) = \prod_{i \in A}(1+(q-1)x_i).
\]
Substituting \(w_i\) yields \(1+(q-1)w_i = q\) if \(u_i = 0\), and \(0\) if \(u_i \neq 0\). Thus,
\[
\mathcal{H}_{\Delta_{\mathcal{I}}}(w_1,\ldots,w_n) = q^{|A|}\delta(u|A).
\]
Substituting \(\mathcal{H}_{\Delta_{\mathcal{I}}}(w_1,\ldots,w_n) = q^{|A|}\delta(u|A)\) into the general weight formula \eqref{eq:weight-Cf}, we obtain
\[
\wt(c_f(s,u)) = (q-1)q^{n-1}(1 - \delta_{0,u}) + (1 - \delta_{s,0})\bigl( q^{|A|}\delta(u|A) - 1 \bigr).
\]
(2) Suppose that \(B \neq \emptyset\). By Lemma \ref{length}, \(|\Delta_{\mathcal{I}}| = q^m + (q-1)^m(q^{|B|} - 1)\). The generating function is
\[
\mathcal{H}_{\Delta_{\mathcal{I}}}(x_1,\ldots,x_n) = \prod_{i \in [m]}(1+(q-1)x_i) + (q-1)^m \prod_{i \in [m]}x_i \left( \prod_{j \in B}(1+(q-1)x_j) - 1 \right).
\]
Let \(u=(v,w)\) with \(v \in \mathbb{F}_q^m\), \(w \in \mathbb{F}_q^{n-m}\), and let \(r = \wt(v)\). Substituting \(w_i\) gives
\[
\mathcal{H}_{\Delta_{\mathcal{I}}}(w_1,\ldots,w_n) = q^m\delta(u|[m]) + (-1)^r(q-1)^{m-r}\left(q^{|B|}\delta(u|B) - 1\right).
\]
We distinguish the following cases based on the values of \(s\), \(v\), and \(\delta(u|B)\):
\begin{enumerate}[label=\textbf{Case \arabic*:}, leftmargin=*]

    \item \(s=0\) and \(u=0\). The weight is \(0\).

    \item \(s=0\) and \(u\neq 0\). The weight is \((q-1)q^{n-1}\).

    \item \(s\neq 0\) and \(u=0\). 
 The weight is \(|\Delta_{\mathcal{I}}|-1 = q^m+(q-1)^m(q^{|B|}-1)-1\).

 \item \(s\neq 0\), \(u\neq 0\), and \(v=0\) (i.e., \(r=0\), \(\delta(u|[m])=1\)).
\begin{enumerate}[label=-, leftmargin=0pt]
\item If \(\delta(u|B)=1\), then \(\mathcal{H} = |\Delta_{\mathcal{I}}|\). The weight is \((q-1)q^{n-1} + q^m + (q-1)^m(q^{|B|}-1) - 1\). 

\item If \(\delta(u|B)=0\), then \(\mathcal{H} = q^m - (q-1)^m\). The weight is \((q-1)q^{n-1} + q^m - (q-1)^m - 1\). The number of such \(u\) is \((q^{|B|}-1)q^{n-m-|B|}\). Multiplying by \(q-1\) gives the frequency \((q-1)(q^{|B|}-1)q^{n-m-|B|}\).
 \end{enumerate}

    \item \(s\neq 0\), \(u\neq 0\), and \(v\neq 0\) (i.e., \(r \geq 1\), \(\delta(u|[m])=0\)).
\begin{enumerate}[label=-, leftmargin=0pt]
        \item If \(\delta(u|B)=1\), then \(\mathcal{H} = (-1)^r(q-1)^{m-r}(q^{|B|}-1)\). The weight is \((q-1)q^{n-1} + (-1)^r(q-1)^{m-r}(q^{|B|}-1) - 1\). 

        \item If \(\delta(u|B)=0\), then \(\mathcal{H} = -(-1)^r(q-1)^{m-r}\). The weight is \((q-1)q^{n-1} - (-1)^r(q-1)^{m-r} - 1\). 
\end{enumerate}

\end{enumerate}
\end{proof}

\begin{table}[htbp]

\centering

\caption{Weight distribution in Theorem~\ref{T3}(2): $B\neq\emptyset$}

\label{tab:c_f_case2}

\renewcommand{\arraystretch}{1.5}

\small 

\begin{tabular}{@{}cc@{}}
\toprule
\textbf{Weight} & \textbf{Frequency} \\
\midrule
\(0\) & \(1\) \\
\((q-1)q^{n-1}\) & \(q^n-1\) \\
\(q^m+(q-1)^m(q^{|B|}-1)-1\) & \(q-1\) \\
\((q-1)q^{n-1} + q^m + (q-1)^m(q^{|B|}-1) - 1\) & \((q-1)(q^{n-m-|B|}-1)\) \\
\((q-1)q^{n-1} + q^m - (q-1)^m - 1\) & \((q-1)(q^{|B|}-1)q^{n-m-|B|}\) \\
\((q-1)q^{n-1} + (-1)^r(q-1)^{m-r}(q^{|B|}-1) - 1\) & \(\binom{m}{r}(q-1)^{r+1}q^{n-m-|B|}\) \\
\((q-1)q^{n-1} - (-1)^r(q-1)^{m-r} - 1\) & \(\binom{m}{r}(q-1)^{r+1}(q^{|B|}-1)q^{n-m-|B|}\) \\
\bottomrule
\end{tabular}

\smallskip

\begin{minipage}{\textwidth}

\footnotesize
Here \(r=\wt(v)\) ranges over \(1 \le r \le m\). Each value of \(r\) gives a separate entry.  If different entries yield the same weight, their frequencies are added together.

\end{minipage}

\end{table}

\begin{remark}
Since $q>2$, every nonzero indicator function taking values in $\{0,1\}$ is nonlinear. Thus, whenever $f\neq0$, the evaluation map $(s,u)\mapsto c_f(s,u)$ is injective. Since the order ideal $I$ in Theorem~\ref{T3} is nonempty, we have $f\neq0$ in all cases.
Consequently,
\(
\dim(\mathcal C_f)=n+1.
\)
The length of $\mathcal C_f$ is $q^n-1$.
\end{remark}

\begin{theorem}\label{T4}

Consider a two-level hierarchical poset \(\mathbb{H}(m,n)\).  Let \(I_1 = A_1 \cup B_1\) and \(I_2 = A_2 \cup B_2\) be two distinct order ideals of \(\mathbb{H}(m,n)\), where \(A_i \subseteq [m]\), \(B_i \subseteq [n] \setminus [m]\), \(i = 1, 2\) and \(I_1 \nsubseteq I_2\), \(I_2 \nsubseteq I_1\). Set \(\mathcal{I} = \{I_1, I_2\}\).

\begin{enumerate}

    \item[(1)] If \(B_1=B_2 = \emptyset\), then the weight distribution of the code \(\mathcal C_f\) is given in Table~\ref{tab:q_dual_1}.

    \item[(2)] If \(B_1 \neq \emptyset\) and \(B_2 \neq \emptyset\), then the weight distribution of the code \(\mathcal C_f\) is given in Table~\ref{tab:q_dual_2}.

\end{enumerate}

\end{theorem}

\begin{proof}
(1) Suppose  \(B_1=B_2=\emptyset\).
In this case, \(I_1=A_1\) and \(I_2=A_2\). By Lemma \ref{rong}, we have
\[
|\Delta_{\mathcal{I}}| = q^{|A_1|}+q^{|A_2|} - q^{|A_1\cap A_2|}.
\]
Moreover, the generating function is
\[
\begin{aligned}
\mathcal{H}_{\Delta_{\mathcal{I}}}(w_1,\ldots,w_n) ={}& q^{|A_1|}\delta(u|A_1) + q^{|A_2|}\delta(u|A_2) - q^{|A_1\cap A_2|}\delta(u|A_1\cap A_2).
\end{aligned}
\]
The cases \(s=0,u=0\), \(s=0,u\neq0\), and \(s\neq0,u=0\) give the first three entries of Table~\ref{tab:q_dual_1}. For \(s\neq0\) and \(u\neq0\), put \(C=A_1\cap A_2\) and \(U=A_1\cup A_2\). We distinguish the following five cases.

\begin{enumerate}[label=\textbf{Case \arabic*:}, leftmargin=*]

\item If \(\delta(u|A_1)=\delta(u|A_2)=1\), then \(\mathcal{H}_{\Delta_{\mathcal{I}}} = q^{|A_1|}+q^{|A_2|}-q^{|C|}\). 

\item If \(\delta(u|A_1)=1\) and \(\delta(u|A_2)=0\), then \(\delta(u|C)=1\), and \(\mathcal{H}_{\Delta_{\mathcal{I}}} = q^{|A_1|}-q^{|C|}\). 

\item If \(\delta(u|A_1)=0\) and \(\delta(u|A_2)=1\), then \(\delta(u|C)=1\), and \(\mathcal{H}_{\Delta_{\mathcal{I}}} = q^{|A_2|}-q^{|C|}\). 

\item If \(\delta(u|A_1)=\delta(u|A_2)=0\) and \(\delta(u|C)=1\), then \(\mathcal{H}_{\Delta_{\mathcal{I}}} = -q^{|C|}\). 

\item If \(\delta(u|C)=0\), then \(\mathcal{H}_{\Delta_{\mathcal{I}}}=0\).

\end{enumerate}

\begin{table}[htbp]

\centering

\caption{Weight distribution of \(\mathcal{C}_f\) for \(B_1=B_2=\emptyset\) over \(\mathbb{F}_q\)}

\label{tab:q_dual_1}

\renewcommand{\arraystretch}{1.5}

\small

\begin{tabular}{cc}
\toprule
\textbf{Weight} & \textbf{Frequency} \\
\midrule
\(0\) & \(1\) \\
\((q-1)q^{n-1}\) & \(q^n-1\) \\
\(q^{|A_1|}+q^{|A_2|}-q^{|A_1\cap A_2|}-1\) & \(q-1\) \\
\((q-1)q^{n-1}+q^{|A_1|}+q^{|A_2|}-q^{|A_1\cap A_2|}-1\) & \((q-1)(q^{n-|A_1\cup A_2|}-1)\) \\
\((q-1)q^{n-1}+q^{|A_1|}-q^{|A_1\cap A_2|}-1\) & \((q-1)q^{n-|A_1\cup A_2|}(q^{|A_2\setminus A_1|}-1)\) \\
\((q-1)q^{n-1}+q^{|A_2|}-q^{|A_1\cap A_2|}-1\) & \((q-1)q^{n-|A_1\cup A_2|}(q^{|A_1\setminus A_2|}-1)\) \\
\((q-1)q^{n-1}-q^{|A_1\cap A_2|}-1\) & \((q-1)q^{n-|A_1\cup A_2|}(q^{|A_1\setminus A_2|}-1)(q^{|A_2\setminus A_1|}-1)\) \\
\((q-1)q^{n-1}-1\) & \((q-1)(q^{|A_1\cap A_2|}-1)q^{n-|A_1\cap A_2|}\) \\
\bottomrule
\end{tabular}

\end{table}

(2) Suppose \(B_1\neq\emptyset\) and \(B_2\neq\emptyset\).
Since \(I_1\) and \(I_2\) are order ideals of \(\mathbb{H}(m,n)\), we have \(A_1=A_2=[m]\). Consequently, \(I_1\cap I_2=[m]\cup(B_1\cap B_2)\).
Put \(T=q^{|B_1|}+q^{|B_2|}-q^{|B_1\cap B_2|}-1\). By Lemma \ref{rong}, we obtain \(|\Delta_{\mathcal{I}}| = q^m+(q-1)^mT\).
For convenience, write
\[
\Delta_B(u) = q^{|B_1|}\delta(u|B_1) + q^{|B_2|}\delta(u|B_2) - q^{|B_1\cap B_2|}\delta(u|B_1\cap B_2) - 1.
\]
Let \(u=(v,w)\), where \(v=(u_1,\ldots,u_m)\) and \(w=(u_{m+1},\ldots,u_n)\), and put \(r=\wt(v)\). The generating function becomes
\[
\mathcal{H}_{\Delta_{\mathcal{I}}}(w_1,\ldots,w_n) = q^m\delta(u|[m]) + (-1)^r(q-1)^{m-r}\Delta_B(u).
\]
The cases \(s=0,u=0\), \(s=0,u\neq0\), and \(s\neq0,u=0\) give the first three entries of Table~\ref{tab:q_dual_2}.
For \(s\neq0\) and \(u\neq0\), we distinguish according to the value of \(r\).
\begin{enumerate}[label=\textbf{Case \arabic*:}, leftmargin=*]

\item Suppose \(r=0\). Then \(\delta(u|[m])=1\) and \((-1)^r(q-1)^{m-r}=(q-1)^m\). There are five possible cases for the restrictions of \(u\) to \(B_1\) and \(B_2\).
\begin{enumerate}[label=-, leftmargin=0pt]

    \item If \(\delta(u|B_1)=\delta(u|B_2)=1\), then \(\Delta_B(u)=T\). 

    \item If \(\delta(u|B_1)=1\) and \(\delta(u|B_2)=0\), then \(\Delta_B(u)=q^{|B_1|}-q^{|B_1\cap B_2|}-1\). 

    \item If \(\delta(u|B_1)=0\) and \(\delta(u|B_2)=1\), then \(\Delta_B(u)=q^{|B_2|}-q^{|B_1\cap B_2|}-1\). 

    \item If \(\delta(u|B_1)=\delta(u|B_2)=0\) and \(\delta(u|B_1\cap B_2)=1\), then \(\Delta_B(u)=-q^{|B_1\cap B_2|}-1\). 

    \item If \(\delta(u|B_1\cap B_2)=0\), then \(\Delta_B(u)=-1\).

\end{enumerate}

\item Suppose \(r\geq1\). Then \(\delta(u|[m])=0\). For each fixed \(r\), there are \(\binom{m}{r}(q-1)^r\) possible vectors \(v\) with \(\wt(v)=r\). Together with the \(q-1\) choices of \(s\neq0\), the common factor in the frequencies is \(\binom{m}{r}(q-1)^{r+1}\). The same five cases for \(\Delta_B(u)\) as above apply.
\end{enumerate}
Combining all the cases yields Table~\ref{tab:q_dual_2}.

\begin{table}[htbp]

\centering

\caption{Weight distribution of \(\mathcal{C}_f\) for \(B_1\neq\emptyset\) and \(B_2\neq\emptyset\) over \(\mathbb{F}_q\)}

\label{tab:q_dual_2}

\renewcommand{\arraystretch}{1.5}

\footnotesize

\begin{tabular}{cc}
\toprule
\textbf{Weight} & \textbf{Frequency} \\
\midrule
\(0\) & \(1\) \\
\((q-1)q^{n-1}\) & \(q^n-1\) \\
\(q^m+(q-1)^mT-1\) & \(q-1\) \\
\((q-1)q^{n-1}+q^m+(q-1)^mT-1\) & \((q-1)(q^{n-m-|B_1\cup B_2|}-1)\) \\
\((q-1)q^{n-1}+q^m+(q-1)^m(q^{|B_1|}-q^{|B_1\cap B_2|}-1)-1\) & \((q-1)q^{n-m-|B_1\cup B_2|}(q^{|B_2\setminus B_1|}-1)\) \\
\((q-1)q^{n-1}+q^m+(q-1)^m(q^{|B_2|}-q^{|B_1\cap B_2|}-1)-1\) & \((q-1)q^{n-m-|B_1\cup B_2|}(q^{|B_1\setminus B_2|}-1)\) \\
\((q-1)q^{n-1}+q^m-(q-1)^m(q^{|B_1\cap B_2|}+1)-1\) & \((q-1)q^{n-m-|B_1\cup B_2|}(q^{|B_1\setminus B_2|}-1)(q^{|B_2\setminus B_1|}-1)\) \\
\((q-1)q^{n-1}+q^m-(q-1)^m-1\) & \((q-1)(q^{|B_1\cap B_2|}-1)q^{n-m-|B_1\cap B_2|}\) \\
\((q-1)q^{n-1}+(-1)^r(q-1)^{m-r}T-1\) & \(\binom{m}{r}(q-1)^{r+1}q^{n-m-|B_1\cup B_2|}\) \\
\((q-1)q^{n-1}+(-1)^r(q-1)^{m-r}(q^{|B_1|}-q^{|B_1\cap B_2|}-1)-1\) & \(\binom{m}{r}(q-1)^{r+1}q^{n-m-|B_1\cup B_2|}(q^{|B_2\setminus B_1|}-1)\) \\
\((q-1)q^{n-1}+(-1)^r(q-1)^{m-r}(q^{|B_2|}-q^{|B_1\cap B_2|}-1)-1\) & \(\binom{m}{r}(q-1)^{r+1}q^{n-m-|B_1\cup B_2|}(q^{|B_1\setminus B_2|}-1)\) \\
\((q-1)q^{n-1}-(-1)^r(q-1)^{m-r}(q^{|B_1\cap B_2|}+1)-1\) & \(\binom{m}{r}(q-1)^{r+1}q^{n-m-|B_1\cup B_2|}(q^{|B_1\setminus B_2|}-1)(q^{|B_2\setminus B_1|}-1)\) \\
\((q-1)q^{n-1}-(-1)^r(q-1)^{m-r}-1\) & \(\binom{m}{r}(q-1)^{r+1}(q^{|B_1\cap B_2|}-1)q^{n-m-|B_1\cap B_2|}\) \\
\bottomrule
\end{tabular}

\vspace{2mm}

\begin{minipage}{\textwidth}

\scriptsize

Here \(T=q^{|B_1|}+q^{|B_2|}-q^{|B_1\cap B_2|}-1\), and \(r=\wt(v)\) ranges over \(1\leq r\leq m\). Each value of \(r\) gives a separate entry. If different entries yield the same weight, their frequencies are added together.

\end{minipage}

\end{table}

\end{proof}

\section{Minimal linear codes}

In this section, based on the weight distribution results established in Section 5, we investigate infinite families of \(q\)-ary linear codes. Specifically, we derive conditions under which these codes are minimal, and analyze their distance optimality for certain parameter configurations over the finite field \(\mathbb{F}_q\).

\begin{theorem}\label{minimal1}
Let \(C_D\) be the \(q\)-ary linear code constructed in Theorem~\ref{T1}(2). If
\(
B=[n]\setminus[m], m\geq 3, |B|\geq 2
\)
and 
\(
q^{m-1}\leq 2(q-1)^{m-1},
\)
then \(C_D\) is a minimal linear code that does not satisfy the Ashikhmin--Barg condition.
\end{theorem}

\begin{proof}
Let
\(
 b=|B|.
\)
Since $B=[n]\setminus[m]$, we have $n=m+b$ and
\[
D=S_m\times(\mathbb F_q^b\setminus\{0\}),
\]
where
\[
S_m=\{\xi\in\mathbb F_q^m:
\xi_i=0\text{ for some }i\in[m]\}.
\]
By Theorem~\ref{T1}(2), we have $\dim(C_D)=n$. Hence the map $u\mapsto c_{D,u}$ is injective. By Lemma~\ref{lem:minimality}, it suffices to show that for any two linearly independent vectors $u,v\in\mathbb F_q^n$,
\[
\sum_{\lambda\in\mathbb F_q^*}
\wt(c_{D,u+\lambda v})
\neq
(q-1)\wt(c_{D,u})-\wt(c_{D,v}).
\]

We first prove that, for every nonzero
$u\in\mathbb F_q^n$,
\begin{equation}\label{eq:span}
\operatorname{span}_{\mathbb F_q}(D\cap\ker u)=\ker u=\{x\in\mathbb F_q^n\ |\ u\cdot x=0\}.
\end{equation}
Write $u=(\alpha,\beta)$, where $\alpha\in\mathbb F_q^m$ and $\beta\in\mathbb F_q^b$.
If $\beta\neq0$, then the projection of
$D\cap\ker u$ onto the first $m$ coordinates is $S_m$, since $b\geq2$ guarantees that
every equation
\(
\beta\cdot z=-\alpha\cdot\xi
\)
has a nonzero solution. Moreover,
\[
\{0\}\times(\ker\beta\setminus\{0\})
\subseteq D\cap\ker u.
\]
Since $S_m$ spans $\mathbb F_q^m$, these facts
imply that
\[
\dim\operatorname{span}_{\mathbb F_q}( D\cap\ker u) \geq m+b-1=\dim\ker u.
\]
Thus \eqref{eq:span} holds. If $\beta=0$, then $\alpha\neq0$ and
\[
D\cap\ker u
=(S_m\cap\ker\alpha)
\times(\mathbb F_q^b\setminus\{0\}).
\]
Choose $j$ such that $\alpha_j\neq0$. The vectors
\[
e_i-\frac{\alpha_i}{\alpha_j}e_j,
\qquad i\neq j,
\]
belong to $S_m\cap\ker\alpha$ and span
$\ker\alpha$, since $m\geq3$.
Also, $\mathbb F_q^b\setminus\{0\}$ spans
$\mathbb F_q^b$.
Therefore,
\[
\operatorname{span}_{\mathbb F_q}(D\cap\ker u)=\ker\alpha\times\mathbb F_q^b=\ker u.
\]

Now let $u,v$ be linearly independent. By the weight-counting identity,
\[
\begin{aligned}
&\sum_{\lambda\in\mathbb F_q^*}
\wt(c_{D,u+\lambda v})
-(q-1)\wt(c_{D,u})+\wt(c_{D,v})=q\big|\{x\in D:
u\cdot x=0,\ v\cdot x\neq0\}\big|.
\end{aligned}
\]
The set on the right is nonempty. Otherwise, $v$ would vanish on $D\cap\ker u$, and hence, by \eqref{eq:span}, on $\ker u$.
This would imply that $v$ is a scalar multiple of $u$, a contradiction. Consequently,
\[
\sum_{\lambda\in\mathbb F_q^*}
\wt(c_{D,u+\lambda v})
>
(q-1)\wt(c_{D,u})-\wt(c_{D,v}).
\]
By Lemma~\ref{lem:minimality}, $C_D$ is minimal.

The length of the code is
\[
N=(q^b-1)(q^m-(q-1)^m).
\]
By the weight distribution in Theorem~\ref{T1}(2), the minimum nonzero weight is
\[
w_{\min}
=\frac{q-1}{q}
\left[N-(q-1)^{m-1}(q^b-1)\right],
\]
while the maximum weight is
\[
w_{\max}
=\frac{q-1}{q}
\left[N+(q-1)^{m-2}(q^b-1)\right].
\]
Indeed, the minimum is attained when $\wt(v)=1$ and $w=0$, whereas the maximum is attained when $\wt(v)=2$ and $w=0$. For the latter, we use
\[
q^m-(q-1)^m
\leq(q+1)(q-1)^{m-1}
\leq(q-1)^{m-2}(q^b-1),
\]
which follows from $b\geq2$ and
$q^{m-1}\leq2(q-1)^{m-1}$. 
It follows that
\[
\begin{aligned}
\frac{w_{\min}}{w_{\max}}
-\frac{q-1}{q}
&=
\frac{N-(q-1)^{m-1}(q^b-1)}
{N+(q-1)^{m-2}(q^b-1)}-\frac{q-1}{q}=
\frac{(q^b-1)\big(q^{m-1}-2(q-1)^{m-1}\big)}
{N+(q-1)^{m-2}(q^b-1)}\leq0.
\end{aligned}
\]
Therefore,
\(
\frac{w_{\min}}{w_{\max}}
\leq\frac{q-1}{q},
\)
and $ C_D$ does not satisfy the Ashikhmin--Barg condition.
\end{proof}

\begin{theorem}\label{minimal2}
Let $q>2$ and let $m\geq 2$. Under the assumptions of Theorem~\ref{T2}(2), the code $C_D$ satisfies the $q$-ary Ashikhmin--Barg condition and hence is a minimal linear code.
\end{theorem}

\begin{proof}
Since $B_1,B_2\neq\emptyset$, we have $A_1=A_2=[m]$. Set
\[
h=|B_1\cap B_2|,\ 
a=|B_1\setminus B_2|,\ 
b=|B_2\setminus B_1|,\ 
c=n-m-|B_1\cup B_2|.
\]
The incomparability of $I_1$ and $I_2$
implies $a,b\geq1$. Let
\[
T=q^{h+a}+q^{h+b}-q^h-1.
\]
Then the length of $C_D$ is
\[
N=q^n-q^m-(q-1)^mT.
\]
By the weight distribution, the minimum nonzero weight is
\[
w_{\min}
=\frac{q-1}{q}
\left(N-(q-1)^{m-1}T\right).
\]

We distinguish two cases.

\medskip
\noindent
\textbf{Case 1: $c\geq1$.}
In this case,
\[
w_{\max}=(q-1)q^{n-1}.
\]
Consequently,
\[
\begin{aligned}
\frac{w_{\min}}{w_{\max}}
-\frac{q-1}{q}
&=\frac{q^{n-1}-q^m
-q(q-1)^{m-1}T}{q^n}.
\end{aligned}
\]
Since $n=m+h+a+b+c$ and $c\geq1$,
we have
\[
q^{n-1}\geq q^{m+h+a+b}.
\]
Moreover, $q(q-1)^{m-1}\leq q^m$.
Thus
\[
\begin{aligned}
&q^{n-1}-q^m-q(q-1)^{m-1}T\geq q^m(q^{h+a+b}-1-T)=q^{m+h}(q^a-1)(q^b-1)>0.
\end{aligned}
\]
Therefore,
\[
\frac{w_{\min}}{w_{\max}}>
\frac{q-1}{q}.
\]

\medskip
\noindent
\textbf{Case 2: $c=0$ and $m\geq2$.}
Without loss of generality, assume $a\leq b$. Since $q\geq3$, we have
\[
\begin{aligned}
4(q^{h+b}-q^h-1)-T&=3q^{h+b}-q^{h+a}-3q^h-3\\
\geq &2q^{h+b}-3q^h-3\\
\geq &(2q-3)q^h-3\\
\geq &0.
\end{aligned}
\]
Hence
\[
T\leq (q-1)^2(q^{h+b}-q^h-1).
\]
It follows from the weight distribution that
\[
w_{\max}
=\frac{q-1}{q}
\left[
N+q^m+(q-1)^m
(q^{h+b}-q^h-1)
\right].
\]
Therefore, the inequality
\[
\frac{w_{\min}}{w_{\max}}>
\frac{q-1}{q}
\]
is equivalent to
\[
q^{m+h+a+b}-q^{m+1}>(q-1)^{m-1}
\left[
(2q-1)T
+(q-1)^2(q^{h+b}-q^h-1)
\right].
\]
Since $m\geq2$ and
\(
\left(\frac{q-1}{q}\right)^{m-2}
\leq1,
\)
it is  sufficient to prove that
\[
q^2(q^{h+a+b}-q)>(q-1)
\left[
(2q-1)T
+(q-1)^2(q^{h+b}-q^h-1)
\right].
\]
Indeed, substituting the expression for $T$, the difference between the two sides of the inequality becomes 
\[
\begin{aligned}
&q^2(q^{h+a+b}-q)-(q-1)
\left[
(2q-1)T
+(q-1)^2(q^{h+b}-q^h-1)
\right]\\
=&q^h\Big[
q(2q-1)+q^2(q^b-q) +(q^a-q)
\big(q^{b+2}-(q-1)(2q-1)\big)
\Big]-q^2\\
\geq &q(2q-1)-q^2\\
=&q(q-1)>0,
\end{aligned}
\]
where we used $a,b\geq1$ and
\[
q^{b+2}>(q-1)(2q-1).
\]
Thus the inequality holds. Hence,
\[
\frac{w_{\min}}{w_{\max}}>
\frac{q-1}{q}.
\]

In both cases, $C_D$ satisfies the $q$-ary Ashikhmin--Barg condition. Therefore, $C_D$ is a minimal linear code.
\end{proof}

\begin{theorem}\label{minimal3}
Let $\mathcal{C}_f$ be the $q$-ary linear code constructed from Theorem~\ref{T3}(2), where $q>2$.
If $m\ge3$, 
\(
B=[n]\setminus[m],
\)
and
\(
|B|\ge2,
\)
then $\mathcal{C}_f$ is a minimal linear code which violates the Ashikhmin--Barg condition.
\end{theorem}

\begin{proof}
Write
\[
x=(\xi,z)\in\mathbb{F}_q^m\times\mathbb{F}_q^b
\]
where \(b=|B|\).  Since $I=[n]$, we have
\[
D=\{(\xi,z):z\neq0,\ 
\xi_i=0\text{ for at least one }i\in[m]\}.
\tag{6.2}
\]
In the following, we divide the proof into five steps.

\medskip
\noindent
\textbf{Step 1. A spanning property.}

We claim that for \(u\neq0\), 
\(
\langle D\cap\ker u\rangle=\ker u,
\)
Let
\(
u=(\alpha,\beta),
\)
where
\(
\alpha\in\mathbb{F}_q^m, 
\beta\in\mathbb{F}_q^b.
\)

\smallskip
\noindent
\textit{\bf{Case 1:} $\beta=0$.} In this case, $\alpha\neq0$. Choose $j\in[m]$ such that
$\alpha_j\neq0$. For every $i\neq j$, put
\[
v_i=e_i-\frac{\alpha_i}{\alpha_j}e_j.
\]
Then 
\(v_i\in \ker\alpha\). Moreover, since $m\geq3$, each $v_i$ has at least one
zero coordinate.
Fix any $z_0\in\mathbb{F}_q^b\setminus\{0\}$. By (6.2), we have
\(
(v_i,z_0),(0,z_0)\in D\cap\ker u.
\)
Therefore
\(
(v_i,0)\in\langle D\cap\ker u\rangle.
\)
Since the vectors
\(\{v_i:\ i\neq j\}\)
form a basis of $\ker\alpha$, we obtain
\[
\ker\alpha\times\{0\}
\subseteq\operatorname{span}_{\mathbb F_q}(D\cap\ker u).
\]
On the other hand,
\(
\{0\}\times
(\mathbb{F}_q^b\setminus\{0\})
\subseteq D\cap\ker u.
\)
Hence
\[
\{0\}\times\mathbb{F}_q^b
\subseteq\operatorname{span}_{\mathbb F_q}(D\cap\ker u).
\]
Combining these inclusions, we get
\[
\begin{aligned}
\operatorname{span}_{\mathbb F_q}(D\cap\ker u)
&=\ker\alpha\times\mathbb{F}_q^b=\ker u.
\end{aligned}
\]

\smallskip
\noindent
\textit{\bf{Case 2:} $\beta\neq0$.} Since $b\geq2$, there exists
\(
h\in\ker\beta\setminus\{0\}.
\)
For $i=1,2$, define
\[
V_i=
\left\{
(\xi,z):
\xi_i=0,\ 
\alpha\cdot\xi+\beta\cdot z=0
\right\}.
\]

We first show that
\(
V_i\subseteq\operatorname{span}_{\mathbb F_q}(D\cap\ker u).
\) Indeed, let $(\xi,z)\in V_i$. If $z\neq0$, then $(\xi,z)\in D\cap\ker u$. If $z=0$, then
\[
(\xi,0)=(\xi,h)+(0,-h).
\]
Since $h\neq0$ and $\beta\cdot h=0$, both
\(
(\xi,h), (0,-h)
\)
belong to $D\cap\ker u$. Thus
\(
V_i\subseteq \operatorname{span}_{\mathbb F_q}(D\cap\ker u).
\)

We next prove that
\(
V_1+V_2=\ker u.
\)
Let $(\xi,z)\in\ker u$. Write
\[
\xi=\xi^{(1)}+\xi^{(2)},
\]
where
\(
\xi^{(1)}=\xi_2e_2,
\xi^{(2)}=\xi-\xi_2e_2.
\)
Then
\(
\xi_1^{(1)}=0,\  \xi_2^{(2)}=0.
\)
Since $\beta\neq0$,  there exists $z_1^0$ such that
\(
\beta\cdot z_1^0
=-\alpha\cdot\xi^{(1)}.
\)
Set
\[
z_2^0=z-z_1^0.
\]
Since $(\xi,z)\in\ker u$, we also have
\(
\beta\cdot z_2^0
=-\alpha\cdot\xi^{(2)}.
\) Set
\[
z_1=z_1^0+th,\qquad
z_2=z_2^0-th,
\]
where $t\in\mathbb{F}_q$.
Because $h\neq0$, each of the equations
\[
z_1^0+th=0,\qquad z_2^0-th=0
\]
has at most one solution for $t$.
As $q>2$, we can choose $t$ such that
\(
z_1\neq0, z_2\neq0.
\)
It follows that
\[
(\xi^{(1)},z_1)\in V_1,\qquad
(\xi^{(2)},z_2)\in V_2.
\]
Therefore, 
\[
(\xi,z)
=(\xi^{(1)},z_1)+(\xi^{(2)},z_2)
\in V_1+V_2.
\]
Thus
\(
\ker u\subseteq V_1+V_2.
\)
The reverse inclusion is immediate, so
\(
V_1+V_2=\ker u
\)
and hence
\[
\ker u
=\operatorname{span}_{\mathbb F_q}(V_1+V_2) \subseteq\operatorname{span}_{\mathbb F_q}(D\cap\ker u).
\]
The reverse inclusion is obvious. Therefore,
\[
\operatorname{span}_{\mathbb F_q}(D\cap\ker u)=\ker u.
\]

\noindent
\textbf{Step 2. The set $D$ spans the whole space.}

We also need to prove
\(
\operatorname{span}_{\mathbb F_q}(D)=\mathbb{F}_q^n.
\)
Indeed, fix $z_0\neq0$. We have 
\(
(0,z_0)
\)
and 
\((e_i,z_0)\in D,
( i\in[m]).
\)
Therefore,
\[
(e_i,0)=(e_i,z_0)-(0,z_0)
\in\operatorname{span}_{\mathbb F_q}(D).
\]
Moreover,
\[
\{0\}\times
(\mathbb{F}_q^b\setminus\{0\})
\subseteq D.
\]
Hence all standard basis vectors of $\mathbb{F}_q^n$ belong to $\operatorname{span}_{\mathbb F_q}(D)$. Hence 
\(
\operatorname{span}_{\mathbb F_q}(D)=\mathbb{F}_q^n.
\)

\medskip
\noindent
\textbf{Step 3. Minimality of $\mathcal{C}_f$.}

Suppose that two nonzero codewords satisfy
\[
\operatorname{supp}(c_f(s,u))
\subseteq
\operatorname{supp}(c_f(s',u')).
\]
Since $f(x)=0$ for every $x\in D$, the above support inclusion implies
\(
D\cap\ker u'\subseteq\ker u.
\)

First suppose that $u'\neq0$. Since 
\(
\ker u'
=\operatorname{span}_{\mathbb F_q}(D\cap\ker u')
\), 
we obtain
\[
\ker u'
=\operatorname{span}_{\mathbb F_q}(D\cap\ker u') \subseteq\ker u.
\]
Thus there exists $\lambda\in\mathbb{F}_q$ such that
\(
u=\lambda u'.
\)
We now show that
\(
s=\lambda s'.
\) Write $u'=(\alpha,\beta)$. We claim that there exists
\(
y\in\Delta_I\setminus\{0\}
\)
such that
\[
u'\cdot y=-s'.
\tag{6.3}
\]

If $\alpha\neq0$, choose $\xi\in\mathbb{F}_q^m$
such that
\(
\alpha\cdot\xi=-s'.
\)
When $s'=0$, we choose $\xi\neq0$ in $\ker\alpha$. This is possible since $m\geq3$. Then $y=(\xi,0)$ satisfies (6.3).

If $\alpha=0$, then $\beta\neq0$. When $s'=0$, take any nonzero $\xi$ and put $y=(\xi,0)$. When $s'\neq0$, take
\(
\xi=(1,\cdots,1)
\)
and choose $z$ satisfying
\(
\beta\cdot z=-s'.
\)
Then $z\neq0$, and $y=(\xi,z)$ belongs to
$\Delta_I\setminus\{0\}$.
Thus (6.3) holds in all cases.
Since $f(y)=1$, we have
\[
sf(y)+u\cdot y
=s+\lambda u'\cdot y
=s-\lambda s'.
\]
On the other hand, (6.3) gives
\(
s'f(y)+u'\cdot y=0.
\)
Since \(\operatorname{supp}(c_f(s,u))
\subseteq
\operatorname{supp}(c_f(s',u'))\), the first codeword must also vanish at $y$.
Therefore,
\(
s-\lambda s'=0.
\)
Hence
\[
(s,u)=\lambda(s',u').
\]
Since $c_f(s,u)$ is nonzero, we have $\lambda\neq0$.

Now suppose that $u'=0$. Since the second codeword
is nonzero, we must have $s'\neq0$.
For every $x\in D$, the second codeword vanishes.
Since \(\operatorname{supp}(c_f(s,u))
\subseteq
\operatorname{supp}(c_f(s',u'))\), the first codeword also vanishes on $D$.
Thus
\[
u\cdot x=0,\qquad x\in D.
\]
As
 \(
\operatorname{span}_{\mathbb F_q}(D)=\mathbb{F}_q^n.
\), we obtain $u=0$.
Consequently, $s,s'\neq0$. Therefore
\[
(s,u)=\frac{s}{s'}(s',u').
\]

We have proved that any two nonzero codewords satisfying
a support inclusion are scalar multiples of each other.
Therefore, $\mathcal{C}_f$ is a minimal linear code.

\medskip
\noindent
\textbf{Step 4. Determination of the minimum and maximum weights.}

By the weight distribution of $\mathcal C_f$, the nonzero weights are of the following forms:
\[
(q-1)q^{n-1},
\]
\[
q^m+(q-1)^m(q^b-1)-1,
\]
\[
(q-1)q^{n-1}+q^m-(q-1)^m-1,
\]
\[
(q-1)q^{n-1}
+(-1)^r(q-1)^{m-r}(q^b-1)-1,
\]
and
\[
(q-1)q^{n-1}
-(-1)^r(q-1)^{m-r}-1,
\]
where $1\leq r\leq m$.

Among the weights in the fourth expression, the minimum is attained at $r=1$, giving
\[
(q-1)q^{n-1}
-(q-1)^{m-1}(q^b-1)-1.
\]
Similarly, the minimum in the fifth expression is attained at $r=2$, giving
\[
(q-1)q^{n-1}-(q-1)^{m-2}-1.
\]
Since
\[
(q-1)^{m-1}(q^b-1)
>(q-1)^{m-2},
\]
the former is smaller than the latter. Moreover,
\[
\begin{aligned}
&(q-1)q^{n-1}
-(q-1)^{m-1}(q^b-1)-1-\left[q^m+(q-1)^m(q^b-1)-1\right]\\
=&q\Big[
\big((q-1)q^b-q\big)
\big(q^{m-2}-(q-1)^{m-2}\big) -(q-1)^{m-2}
\Big]\\
>&0.
\end{aligned}
\]
Indeed, since $m\geq3$ and $b\geq2$,
\[
q^{m-2}-(q-1)^{m-2}
\geq(q-1)^{m-3}
\]
and
\[
(q-1)q^b-q>q-1.
\]
Furthermore,
\[
q^m-(q-1)^m-1>0.
\]
Consequently, the minimum nonzero weight is
\[
w_{\min}
=q^m+(q-1)^m(q^b-1)-1.
\tag{6.4}
\]

Next, the maximum in the fourth expression is attained at $r=2$, giving
\[
(q-1)q^{n-1}
+(q-1)^{m-2}(q^b-1)-1.
\]
The maximum in the fifth expression is attained at $r=1$, giving
\[
(q-1)q^{n-1}+(q-1)^{m-1}-1.
\]
Since
\[
q^m-(q-1)^m>(q-1)^{m-1},
\]
we have
\[
(q-1)q^{n-1}+q^m-(q-1)^m-1>
(q-1)q^{n-1}+(q-1)^{m-1}-1.
\]
Therefore,
\[
\begin{aligned}
w_{\max}
={}&(q-1)q^{n-1}+\max\left\{
q^m-(q-1)^m,\,
(q-1)^{m-2}(q^b-1)
\right\}-1.
\end{aligned}
\tag{6.5}
\]
In particular,
\[
w_{\max}>(q-1)q^{n-1}.
\]

\medskip
\noindent
\textbf{Step 5. Violation of the Ashikhmin--Barg condition.}

By (6.4) and (6.5), we have
\[
\begin{aligned}
\frac{w_{\min}}{w_{\max}}
&<
\frac{q^m+(q-1)^m(q^b-1)-1}
{(q-1)q^{n-1}}\\
&=
\frac{q^{1-b}}{q-1}
+\left(\frac{q-1}{q}\right)^{m-1}
(1-q^{-b})-\frac{1}{(q-1)q^{n-1}}\\
&<
\frac{1}{q(q-1)}
+\left(\frac{q-1}{q}\right)^2\\
&<
\frac{q-1}{q}.
\end{aligned}
\]
Here we used $n=m+b$, $m\geq3$,
$b\geq2$, and
\[
(q-1)^2>q,
\]
which holds for $q>2$.

Therefore, $\mathcal C_f$ violates the Ashikhmin--Barg condition. Combining this with Step 3, we conclude that $\mathcal C_f$ is a minimal linear code which violates the Ashikhmin--Barg condition.
\end{proof}

\begin{theorem}\label{minimal4}
Let $\mathcal{C}_f$ be the $q$-ary linear code constructed in Theorem~\ref{T4}(2), where $q>2$. Suppose that
\[
m=2,\ n=4,\ B_1=\{3\},\ B_2=\{4\}.
\]
Then $\mathcal{C}_f$ is a minimal linear code with parameters
\(
[q^4-1,\ 5,\ q^2+2(q-1)^3-1]_q,
\)
which violates the Ashikhmin--Barg condition.
\end{theorem}

\begin{proof}
Write
\[
(\xi,\eta)\in
\mathbb F_q^2\times\mathbb F_q^2,
\]
where $\xi=(\xi_1,\xi_2)$ and
$\eta=(\eta_1,\eta_2)$.
Note that
\[
\Delta_{\mathcal I}
=
\{(\xi,0):\xi\in\mathbb F_q^2\}
\cup
\left\{
(\xi,\eta):
\xi_1\xi_2\neq0,\ 
|\operatorname{supp}(\eta)|=1
\right\}.
\]
Consequently,
\(
|\Delta_{\mathcal I}|
=q^2+2(q-1)^3.
\)
Moreover, 
\[
D=\{(\xi,\eta):\eta_1\eta_2\neq0\text{ or }(\eta\neq0,\ \xi_1\xi_2=0)\}.
\]

\medskip
\noindent
\textbf{Step 1. A spanning property.} 

We claim that  for every nonzero $(p,r)\in\mathbb F_q^2
\times\mathbb F_q^2$,
\[
\operatorname{span}_{\mathbb F_q}(D\cap\ker(p,r))
=\ker(p,r)=
\left\{(\xi,\eta):
p\cdot\xi+r\cdot\eta=0
\right\}.
\tag{6.6}
\]

First, we show that the evaluation map
\[
(s,u)\longmapsto c_f(s,u)
\]
is injective. Suppose that
\[
sf(x)+u\cdot x=0
\]
for every nonzero $x\in\mathbb F_q^4$. 

If $s=0$, then $u=0$. If $s\neq0$ and $u\neq0$, the linear form $u\cdot x$ takes every value in $\mathbb F_q$ on nonzero vectors, whereas $f(x)$ only takes the values $0$ and $1$. This is impossible because $q>2$.

If $s\neq0$ and $u=0$, then $f$ would vanish on every nonzero vector, which is also impossible. Therefore, the evaluation map is injective, and
\(
\dim(\mathcal C_f)=5.
\)

To prove (6.6), let
\(
(a,b)\in
\operatorname{span}_{\mathbb F_q}(D\cap\ker(p,r))^\perp.
\)
Then
\(
a\cdot\xi+b\cdot\eta=0
\)
for every $(\xi,\eta)\in D\cap\ker(p,r)$.

\smallskip
\noindent
\textit{\bf{Case 1:} $r=0$.} In this case, $p\neq0$.
Choose
\(
\gamma\in\mathbb F_q^*\setminus\{1\},
\)
which is possible since $q>2$. Since
\(
(0,(1,1)),\ (0,(1,\gamma))
\in D\cap\ker(p,0),
\)
we obtain
\[
b_1+b_2=0,\qquad
b_1+\gamma b_2=0.
\]
Hence $b=0$. Moreover, for every $\xi\in\ker p$,
\(
(\xi,(1,1))\in D\cap\ker(p,0).
\)
Hence
\(
a\cdot\xi=0.
\)
It follows that $a=\lambda p$ for some
$\lambda\in\mathbb F_q$. Therefore,
\[
(a,b)=\lambda(p,0)=\lambda(p,r).
\]

\smallskip
\noindent
\textit{\bf{Case 2:}} Exactly one coordinate of $r$ is nonzero. Without loss of generality, suppose that
\(
r=(r_1,0),
\)
where
\(
 r_1\neq 0.
\)
If $p=0$, then
\(
(0,(0,1)),
(e_1,(0,1)),
(e_2,(0,1))
\)
all belong to $D\cap\ker(0,r)$. Therefore,
\(
b_2= a_1= a_2=0.
\)
Thus there exists $\lambda\in\mathbb F_q$ such that
\[
(a,b)=\lambda(0,r)=\lambda(p,r).
\]

Now suppose that $p\neq0$. Since
\(
(0,(0,1))\in D\cap\ker(p,r),
\)
we have $b_2=0$. For every $\xi\in\mathbb F_q^2$ satisfying $p\cdot\xi\neq0$, put
\(
\eta=\left(-\frac{p\cdot\xi}{r_1},1\right).
\)
Then $\eta_1\eta_2\neq0$ and
\(
p\cdot\xi+r\cdot\eta=0.
\)
Hence
\(
(\xi,\eta)\in D\cap\ker(p,r).
\)
It follows that
\(
\left(a-\frac{b_1}{r_1}p\right)\cdot\xi=0
\)
for every $\xi$ satisfying $p\cdot\xi\neq0$. Note that the set
\(
\{\xi\in\mathbb F_q^2:p\cdot\xi\neq0\}
\)
spans $\mathbb F_q^2$.  Consequently,
\(
a=\frac{b_1}{r_1}p.
\)
Together with $b_2=0$, this gives
\[
(a,b)=\frac{b_1}{r_1}(p,r).
\]

\smallskip
\noindent
\textit{\bf{Case 3:}} Both coordinates of $r$ are nonzero.
Suppose that $r_1r_2\neq0$. Let 
\(
\eta_0=(r_2,-r_1).
\)
Then
\(
 r\cdot\eta_0=0.
\)
For every $\xi\in\ker p$, we have
\(
(\xi,\eta_0)\in D\cap\ker(p,r).
\)
Taking $\xi=0$ first, we obtain
\(
b\cdot\eta_0=0.
\)
Thus $b=\lambda r$ for some $\lambda\in\mathbb F_q$. Furthermore,
\(
a\cdot\xi=0
\)
for every 
\(
 \xi\in\ker p.
\)
If $p=0$, this implies $a=0$. Hence 
\[
(a,b)=\lambda(p,r).
\]

Now suppose that $p\neq0$. Then there exists $\mu\in\mathbb F_q$ such that
\(
a=\mu p.
\)
Choose $\gamma$ satisfying $p\cdot\gamma\neq0$, and let $c=-p\cdot\gamma$. Since $q>2$, we can choose
\(
\eta_1\in\mathbb F_q^*
\setminus\left\{\frac{c}{r_1}\right\}.
\)
Let
\(
\eta_2=\frac{c-r_1\eta_1}{r_2}.
\)
Then $\eta_1\eta_2\neq0$ and
\(
r\cdot\eta=c=-p\cdot\gamma.
\)
Therefore,
\(
(\gamma,\eta)\in D\cap\ker(p,r).
\)
Using $a=\mu p$ and $b=\lambda r$, we obtain
\[
0=a\cdot\gamma+b\cdot\eta
=\mu p\cdot\gamma+\lambda r\cdot\eta=(\mu-\lambda)p\cdot\gamma.
\]
Since $p\cdot\gamma\neq0$, it follows that
\(
\mu=\lambda.
\)
Thus
\[
(a,b)=\lambda(p,r).
\]

Combining the three cases, we conclude that
\[
\operatorname{span}_{\mathbb F_q}(D\cap\ker(p,r))^\perp
=\langle(p,r)\rangle.
\]
Thus
\(
\dim\left\langle D\cap\ker(p,r)\right\rangle=3.
\)
Since $\ker(p,r)$ also has dimension $3$, we obtain 
\[
\operatorname{span}_{\mathbb F_q}(D\cap\ker(p,r))=\ker(p,r)
\]

\noindent
\textbf{Step 2. The set $D$ spans the whole space.}

We claim that
\[
\operatorname{span}_{\mathbb F_q}(D)=\mathbb F_q^4.
\]
Indeed, choose
$\gamma\in\mathbb F_q^*\setminus\{1\}$.
The four vectors
\[
(0,(1,1)),\quad
(0,(1,\gamma)),\quad
(e_1,(1,1)),\quad
(e_2,(1,1))
\]
belong to $D$ and span $\mathbb F_q^4$.

\medskip
\noindent
\textbf{Step 3. Minimality of $\mathcal{C}_f$.}

Now suppose that two nonzero codewords satisfy
\[
\operatorname{supp}(c_f(s,u))
\subseteq
\operatorname{supp}(c_f(s',u')),
\]
which implies
\(
D\cap\ker u'\subseteq\ker u.
\)
 Write $u'=(p',r')$.

First suppose that $u'\neq0$. Since 
\(
\ker (p',r')
=\operatorname{span}_{\mathbb F_q}(D\cap\ker (p',r'))
\), 
we obtain
\[
\ker u'
=\operatorname{span}_{\mathbb F_q}(D\cap\ker u')
\subseteq\ker u.
\]
Thus there exists $\lambda\in\mathbb{F}_q$ such that
\(
u=\lambda u'.
\)
We now show that
\(
s=\lambda s'.
\)
We claim that there exists
\(
y\in\Delta_{\mathcal I}\setminus\{0\}
\)
such that
\[
u'\cdot y=-s'.
\tag{6.7}
\]

If $p'\neq0$, choose $\xi\in\mathbb F_q^2$ such that
\(
p'\cdot\xi=-s'.
\)
When $s'=0$, choose a nonzero
$\xi\in\ker p'$.
Then $y=(\xi,0)$ satisfies (6.7) and $f(y)=1$. 

If $p'=0$, then $r'\neq0$. When $s'=0$, take
\(
y=((1,1),0).
\)
When $s'\neq0$, choose $j\in\{1,2\}$ such that $r'_j\neq0$, and let
\(
y=\left((1,1),-\frac{s'}{r'_j}e_j\right).
\)
Again, $y$ satisfies (6.7) and $f(y)=1$. Since $f(y)=1$, we have
\[
sf(y)+u\cdot y
=s+\lambda u'\cdot y
=s-\lambda s'.
\]
On the other hand, (6.7) gives
\(
s'f(y)+u'\cdot y=0.
\)
Since \(\operatorname{supp}(c_f(s,u))
\subseteq
\operatorname{supp}(c_f(s',u'))\), the first codeword must also vanish at $y$. Therefore,
\(
s-\lambda s'=0.
\)
Hence
\[
(s,u)=\lambda(s',u').
\]
Since $c_f(s,u)$ is nonzero, we have $\lambda\neq0$.

It remains to consider the case $u'=0$. Since $c_f(s',u')$ is nonzero, we have $s'\neq0$. For every $x\in D$, the second codeword vanishes. Since \(\operatorname{supp}(c_f(s,u))
\subseteq
\operatorname{supp}(c_f(s',u'))\), the first codeword also vanishes on $D$. Since $f(x)=0$ on $D$, we obtain
\[
u\cdot x=0,\qquad x\in D.
\]
As
 \(
\operatorname{span}_{\mathbb F_q}(D)=\mathbb{F}_q^n
\), 
we obtain $u=0$. Consequently, $s,s'\neq0$. Therefore
\[
(s,u)=\frac{s}{s'}(s',u').
\]

We have proved that any two nonzero codewords satisfying a support inclusion are scalar multiples of each other. Therefore, $\mathcal{C}_f$ is a minimal linear code.

\medskip
\noindent
\textbf{Step 4. Determination of the minimum and maximum weights.}

By the weight distribution of $\mathcal C_f$, if $s=0$ and $u\neq0$, then
\[
\wt(c_f(s,u))=(q-1)q^3.
\]
If $s\neq0$ and $u=0$, then
\[
\wt(c_f(s,0))
=q^2+2(q-1)^3-1.
\]

Suppose that $s\neq0$ and $u\neq0$. Multiplying the codeword by $s^{-1}$ does not change its weight, so we may assume $s=1$. Writing $u=(p,r)$, the weight formula gives
\[
\wt(c_f(1,(p,r)))=(q-1)q^3+H(p,r)-1,
\]
where
\[
\begin{aligned}
H(p,r)
={}&q^2\delta_{0,p}
+(-1)^{\wt(p)}(q-1)^{2-\wt(p)}
\bigl(
q\delta_{0,r_1}+q\delta_{0,r_2}-2
\bigr).
\end{aligned}
\]

Considering $\wt(p)=0,1,2$ and the number of zero coordinates of $r$, we obtain
\[
-2(q-1)^2\leq H(p,r) \leq q^2+(q-1)^2(q-2).
\]
The lower bound is attained when $\wt(p)=1$ and $r=0$, while the upper bound is attained when $p=0$ and exactly one coordinate of $r$ is nonzero.

Thus, the minimum weight among the codewords with $s\neq0$ and $u\neq0$ is
\[
(q-1)q^3-2(q-1)^2-1.
\]
A direct calculation gives
\[
(q-1)q^3-2(q-1)^2-1-\bigl[q^2+2(q-1)^3-1\bigr]=q(q-2)(q^2-q+1)>0.
\]
Moreover,
\[
(q-1)q^3-\bigl[q^2+2(q-1)^3-1\bigr]=(q-1)(q^3-2q^2+3q-3)>0.
\]
Therefore, the minimum nonzero weight is
\[
w_{\min}=q^2+2(q-1)^3-1.
\]

Similarly, the maximum nonzero weight is
\[
w_{\max}=(q-1)q^3+q^2+(q-1)^2(q-2)-1=(q-1)(q^3+q^2-2q+3).
\]
Consequently, the parameters of $\mathcal C_f$
are
\[
[q^4-1,\ 5,\ q^2+2(q-1)^3-1]_q.
\]

\medskip
\noindent
\textbf{Step 5. Violation of the Ashikhmin--Barg condition.}

Since
\[
w_{\max}>(q-1)q^3,
\]
we have
\[
\begin{aligned}
\frac{w_{\min}}{w_{\max}}
&<
\frac{q^2+2(q-1)^3-1}{(q-1)q^3}.
\end{aligned}
\]
Furthermore,
\[
\begin{aligned}
&\frac{q-1}{q}
-\frac{q^2+2(q-1)^3-1}{(q-1)q^3}=\frac{(q-1)^3-2}{q^3}>0,
\end{aligned}
\]
because $q>2$. Hence
\[
\frac{w_{\min}}{w_{\max}}
<\frac{q-1}{q}.
\]
Therefore, $\mathcal C_f$ is a minimal linear code which violates the Ashikhmin--Barg condition.
\end{proof}


\begin{thebibliography}{99}

\bibitem{Ashikhmin1998}
A. Ashikhmin, A. Barg,
Minimal vectors in linear codes,
IEEE Trans. Inf. Theory 44 (5) (1998) 2010--2017.

\bibitem{Coulter1997}
R. S. Coulter, R. W. Matthews,
Bent polynomials over finite fields,
Bull. Aust. Math. Soc. 56 (1997) 429--437.

\bibitem{Ding2003}
C. Ding, J. Yuan,
Covering and secret sharing with linear codes,
in: C. S. Calude, M. J. Dinneen, V. Vajnovszki (Eds.),
Discrete Mathematics and Theoretical Computer Science,
Lecture Notes in Computer Science, vol. 2731,
Springer, Berlin, Heidelberg, 2003, pp. 11--25.

\bibitem{Ding2016}
C. Ding,
A construction of binary linear codes from Boolean functions,
Discret. Math. 339 (2016) 2288--2303.

\bibitem{Ding2018}
C. Ding, Z. Heng, Z. Zhou,
Minimal binary linear codes,
IEEE Trans. Inf. Theory 64 (2018) 6536--6545.

\bibitem{Griesmer1960}
J. H. Griesmer,
A bound for error correcting codes,
IBM J. Res. Dev. 4 (5) (1960) 532--542.

\bibitem{Heng2018}
Z. Heng, C. Ding, Z. Zhou,
Minimal linear codes over finite fields,
Finite Fields Appl. 54 (2018) 176--196.

\bibitem{Huffman2003}
W. C. Huffman, V. Pless,
Fundamentals of Error-Correcting Codes,
Cambridge University Press, Cambridge, 2003.

\bibitem{Hyun2020}
J. Y. Hyun, J. Lee, Y. Lee,
Infinite families of optimal linear codes constructed from simplicial complexes,
IEEE Trans. Inf. Theory 66 (11) (2020) 6762--6773.

\bibitem{Hyun2020poset}
J. Y. Hyun, H. K. Kim, Y. Wu, Q. Yue,
Optimal minimal linear codes from posets,
Des. Codes Cryptogr. 88 (2020) 2475--2492.

\bibitem{Lu2021}
W. Lu, X. Wu, X. Cao,
The parameters of minimal linear codes,
Finite Fields Appl. 71 (2021) 101799.

\bibitem{Lu2021poset}
W. Lu, X. Wu, X. Cao,
Minimal linear codes from posets,
Des. Codes Cryptogr. 89 (2021) 1109--1126.

\bibitem{Massey1993}
J. L. Massey,
Minimal codewords and secret sharing,
in: Proc. 6th Joint Swedish-Russian Workshop on Information Theory,
1993, pp. 276--279.

\bibitem{Mesnager2019}
S. Mesnager, F. \\"{O}zbudak, A. S\i nak, G. Cohen,
On $q$-ary plateaued functions over $\mathbb F_q$
and their explicit characterizations,
European J. Combin. 80 (2019) 71--81.

\bibitem{Pelen2026}
R. M. Pelen,
Optimal and minimal $p$-ary linear codes from generalized order ideals of hierarchical posets,
arXiv:2609.33128 (2026).

\bibitem{Shamir1979}
A. Shamir,
How to share a secret,
Commun. ACM 22 (11) (1979) 612--613.

\bibitem{Wu2025}
X. Wu, W. Lu, X. Cao, X. Qin,
Minimal linear codes constructed from hierarchical posets with two levels,
Adv. Math. Commun. 19 (1) (2025) 276--283.

\end{thebibliography}
\end{document}